%% file: MultiShapeManuscript.tex
\documentclass[onefignum,onetabnum]{siamart190516}

\input{ex_shared}
\usepackage{eufrak}

\begin{document}

\maketitle

\begin{abstract}
A numerical framework is developed to solve various types of PDEs on complicated domains, including steady and time-dependent, non-linear and non-local PDEs, with different boundary conditions that can also include non-linear and non-local terms. This numerical framework, called MultiShape, is a class in {\scshape Matlab}, and the software is open source. We demonstrate that MultiShape is compatible with other numerical methods, such as differential--algebraic equation solvers and optimization algorithms. The numerical implementation is designed to be user-friendly, with most of the set-up and computations done automatically by MultiShape and with intuitive operator definition, notation, and user-interface. Validation tests are presented, before we introduce three examples motivated by applications in Dynamic Density Functional Theory and PDE-constrained optimization, illustrating the versatility of the method.
\end{abstract}



\section{Introduction}
Many problems in the applied sciences, including biology, chemical engineering, and physics can be described by (integro-) partial differential equation (PDE) models. These include wide-ranging applications of industrial relevance, including those in drug delivery \cite{PEPPAS201475}, manufacturing \cite[Chapter 2]{jaluria2018mathematical}, and the food industry \cite{foodprocessing1}. Often in these applications, such models need to be solved on complicated domains, to include relevant features of the problem setup, for example in microfluidics \cite{bridle2014application} and milling processes in the pharmaceutical industry \cite{seibert2019milling}. Therefore, efficient numerical methods for solving integro-PDEs on complicated domains are of ever-increasing interest to academic and industrial communities \cite{RaoFEM,Canuto2007,lozinski2003fokker,lozinski2003fast,FAKHARIZADI1,zhu1}.

Perhaps the most popular numerical scheme for such problems is the finite element method (FEM) \cite{BrennerScottFEM}, which is both an active area of applied mathematics research and frequently used in a wide range of applications \cite{RaoFEM}, with several commercial and open-source software solutions available, such as FEniCS \cite{Fenics1,Fenics2} and Abaqus \cite{Abaqus1}.  However, for a large class of PDE models, there is another effective numerical method---the pseudospectral method---which, in particular for modelling problems that involve non-local phenomena, can be more efficient and accurate than the FEM for a comparable computational cost. We are particularly interested in applications in which non-local phenomena enter the model through, for example, convolution terms, on which we focus here, or fractional operators. 
For these models, FEM becomes computationally expensive because one of its key advantages, the frequently obtained sparse matrix structures, is generally lost. Pseudospectral methods, however, do not rely on sparsity and non-local terms can be treated without significant additional computational cost. Moreover, pseudospectral methods exhibit exponential convergence in the number of collocation points for smooth functions and can exploit any smoothness of the PDE solution \cite[Chapter 1]{Boyd1}. For many models that include diffusion and sufficiently `nice' interactions, smooth solutions can be expected, even if such results have not yet been rigorously proven.  Any potential such proofs are expected to be highly challenging.
For example, even with periodic boundary conditions, non-local problems of the form studied here may not have a unique equilibrium, depending on the relative strengths of diffusion and interaction~\cite{carrillo2020long}, which can provide significant numerical challenges in regions of parameter space that lead produce `phase transitions'.

In \cite{GoddardPseudospectralCode1}, a pseudospectral scheme has been introduced to efficiently solve non-local PDEs, including complicated non-local boundary conditions \emph{on simple domains}. The authors provided an open-source software implementation of their method, 2DChebClass, which applies to 1D problems and simple domains in two dimensions, such as rectangles, wedges (sections of an annulus, in polar coordinates), and semi-infinite boxes. Features of this implementation include the straightforward computation of convolution integrals and the user-friendly implementation of complicated, often non-local, boundary conditions, which cannot easily be tackled by standard methods, such as boundary bordering (see e.g., \cite[Chapter 6]{Boyd1} and \cite[Chapter 7]{bibTrefethen}).

A natural extension of pseudospectral methods is the spectral element method (SEM), which can be seen as combining a pseudospectral method with a higher-order finite element method. The basis functions used by SEM, such as Chebyshev or Lagrange polynomials, are typically of higher order than those used in FEM.  Additionally, SEM generally uses Chebyshev--Lobatto collocation grids on each element, as compared to a more standard equispaced grid in FEM. At the intersections between the elements, continuity is enforced by imposing two matching conditions, usually continuity of the solution and its first derivative normal to the intersection, for PDE problems that are second order in space; see \cite[Chapter 22]{Boyd1}. 
There are a number of different, but related, spectral element methods. One can consider the strong form of the PDE and apply matching conditions on the intersections of the elements directly. This is called the patching method, and was first introduced by Orszag \cite{ORSZAG198070}. The Galerkin spectral element method solves the weak form of a PDE model, and was pioneered by Patera \cite{SEMPatera84} using Chebyshev polynomials as basis functions, and later adapted to Lagrange polynomials by Komatitsch and Vilotte \cite{SEMLagrange98}, which is now the standard choice \cite[Chapter 22]{Boyd1}. As shown in \cite{Funaro1}, for elliptic PDEs this approach is essentially equivalent to the patching formulation. A third method for solving a PDE using SEM is based on overlapping elements and was introduced by Morchoisne \cite{Morchoisne1}.

Spectral element methods have been applied to solve a range of problems. Some of the earlier applications were in (classical) fluid dynamics (see \cite[Chapter 13]{Canuto1988} and \cite[Chapter 5]{Canuto2007}) and geosciences \cite{SEMLagrange98,fichtner2009simulation}. Further applications of SEM include fractional PDEs \cite{lischke2019spectral}, relativity \cite{pfeiffer2003multidomain,grandclement2009spectral}, financial modelling \cite{zhu3,zhu2,zhu1}, and optimal control problems with elliptic PDE constraints \cite{chen2017spectral,zhou2011spectral,chen2008legendre}. 
Often SEM is applied as a higher-order finite element method, using a large number of elements with low degree polynomials. However, some work has been done on domain decomposition and `multidomain' methods, e.g., \cite{pfeiffer2003multidomain}, that actively apply SEM to describe domains that are more complicated. Examples of these include modelling blood flow in veins, aircraft engines, or flows around airfoils \cite[Chapter 5]{Canuto2007}, flows around obstacles \cite{lozinski2003fokker}, flows in petroleum reservoirs \cite{TANEJA1}, seismic waves \cite{SEMLagrange98,fichtner2009simulation}, and sound propagation \cite{lin1998least}.
While most existing work does not include non-local integral terms, some advances have been made, see e.g., \cite{lozinski2003fokker,zhu1,FAKHARIZADI1,fakhar2015spectral}. 
These works consider different, mostly one dimensional, integral terms, which were tackled in the context of the (weak form) Galerkin SEM formalism, using either Gaussian quadrature or Fast Fourier Transforms. Additionally, the chosen boundary conditions do not involve convolution integral terms, unlike those for the problems considered here.

The present work develops a multidomain spectral element method, called MultiShape\footnote{We apply the convention that `MultiShape', with capital letters, refers to our methodology as a whole, while a lower-case `multishape' refers to an individual shape (i.e., domain).}, including an open-source software package \cite{MultiShapeCode}, which is able to solve various (integro-)PDE problems on complicated domains. This implementation includes the application of complicated boundary conditions, such as (non-linear, non-local) no-flux conditions, as well as the efficient solution of integro-PDEs, involving the evaluation of convolution integrals. For the problems considered in the present work, this is a significant improvement to the approaches taken in the previous works mentioned above.
The work follows from the pseudospectral code library 2DChebClass \cite{GoddardPseudospectralCode1}, extending this framework to a spectral element method, with minimal additional effort for the user. The flexibility of the method is illustrated by a range of numerical examples in Sections \ref{sec:Validation} and \ref{sec:NumericalResults}, demonstrating one of the main advantages of the implementation: its compatibility with various other numerical methods such as standard ordinary differential equation (ODE) and state-of-the-art optimization algorithms. In line with previous work, see \cite{GoddardPseudospectralCode1,OurPreprint2021}, we present examples that originate from Dynamic Density Functional Theory (DDFT), which describes dynamic properties of systems involving particles suspended in a fluid. Applications of this theory range from mathematical biology, such as in bird flocking or bacteria dynamics \cite{menzel2016dynamical}, through medicine, such as in drug delivery \cite{fang2005kinetics}, to social sciences, such as in opinion dynamics \cite{goddard2022noisy}. 
In DDFT, complicated domains are often modelled by imposing steep potential gradients in parts of the domain, e.g., to model flow through a constriction \cite{zimmermann2016flow}, or around obstacles \cite{rauscher2007dynamic,almenar2011dynamics}, which can be computationally expensive due to large gradients in the system. Here we demonstrate that spectral element methods are a valuable alternative to this method and solutions to DDFT-like problems on complicated domains using SEM are illustrated in Section \ref{sec:NumericalResults}. 

In summary, the present work makes the following significant contributions: 
\begin{itemize}
	\item Open-source software for the computation of a broad range of (integro-)PDE problems on complicated domains.
	\item Straightforward to interface with with other numerical schemes, including differential--algebraic equation (DAE) solvers and optimization algorithms, with problem set-ups being independent of the domain; one may readily re-solve on a range of domains with essentially no additional coding effort.
	\item User-friendly set-up, compatible with the 2DChebClass library and easily adapted to new model problems.
	\item Easy application of naturally-arising complicated, non-linear, non-local boundary conditions, which do not have to be written by the user in a problem-specific way.
	\item Automated constructions of multishapes, whose elements can be defined by a combination of Cartesian and polar coordinates.
	\item Efficient evaluation of convolution integrals on complicated domains, easily applicable to a wide range of kernels.
\end{itemize}
The paper is structured as follows. In Section \ref{sec:NumericalMethods}, the numerical method is introduced in detail. In Section \ref{sec:Validation} the code is validated with some numerical tests. Section \ref{sec:ModelEquations} introduces the three applications considered in detail: an equilibrium model, an integro-PDE model for particle dynamics, and a PDE-constrained optimal control model. Furthermore, the numerical methods that are used in combination with SEM are introduced. Section \ref{sec:NumericalResults} presents some numerical results before we close with some remarks in Section \ref{sec:Conclusions}.

\section{Construction of Numerical Methods}\label{sec:NumericalMethods}

\begin{figure}[h]
	\centering
\includegraphics[scale=0.17]{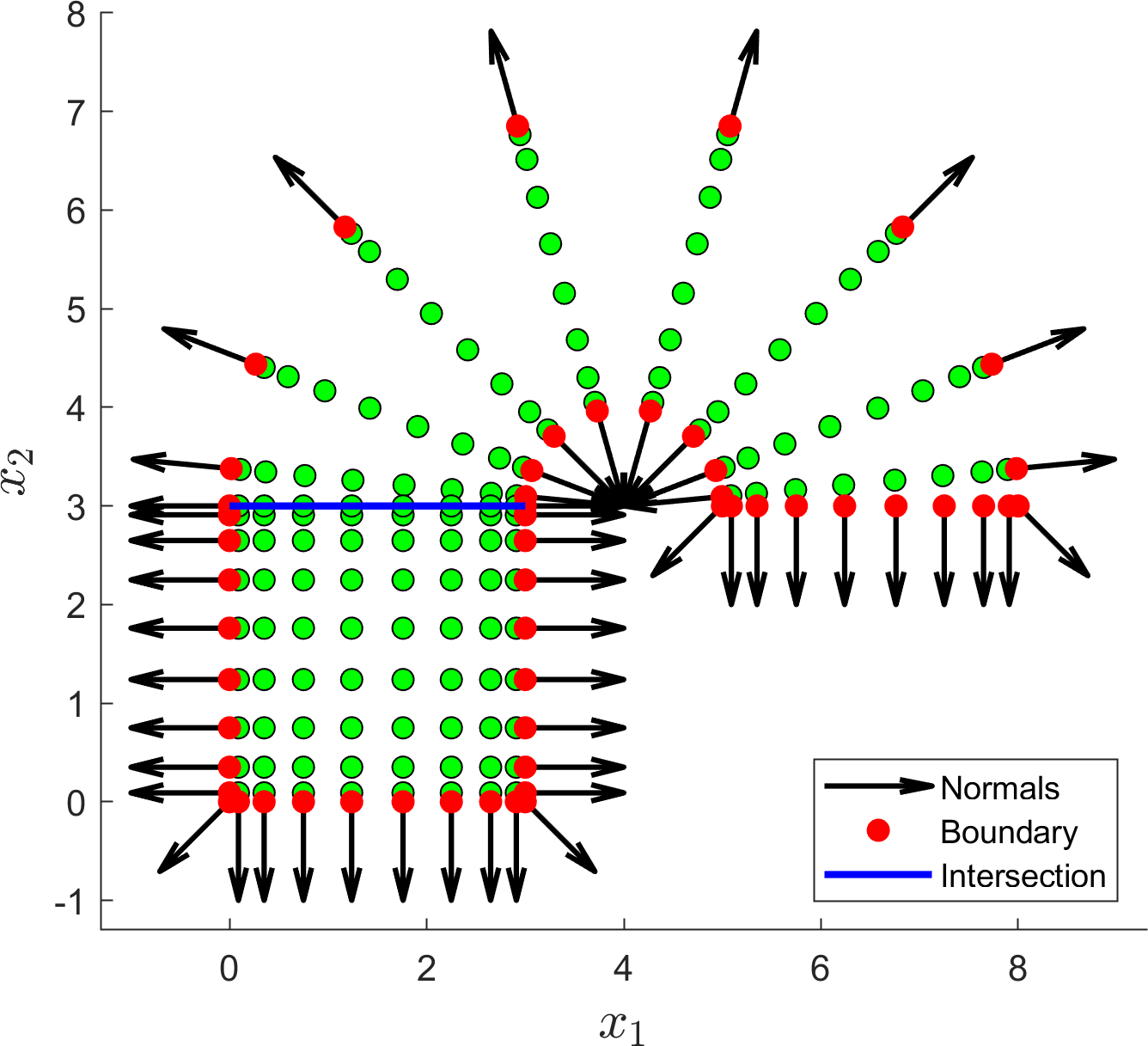}
	\caption{A multishape consisting of a quadrilateral and a wedge element. Displayed are collocation points, highlighted boundary points, the intersection boundary between the two shapes, and outward normal vectors.} 
	\label{fig:MS1}
\end{figure}

In this section a numerical methodology, called MultiShape, is developed to solve various types of PDEs on complicated domains, including steady and time-dependent, non-linear and non-local PDEs, with different boundary conditions that can also include non-linear and non-local terms. Other problems involving differentiation, interpolation, integration, or convolution, such as ODEs or integral equations, can also be tackled using this framework.
MultiShape is a class in {\scshape Matlab}, and is based on 2DChebClass, an object oriented library for solving PDEs on single elements.
A multishape domain is defined by a number $n_e$ of simple elements, denoted by $e_{i}$, $i = 1,...,n_e$. In the present work, quadrilateral and wedge (section of an annulus) elements can be combined to construct the multishape domain (a discretization of the problem-specific physical domain); see Figure \ref{fig:MS1} for an illustration.  This can be straightforwardly extended to other shapes that can be mapped (conformally) from/to the unit square (i.e., those used in standard pseudospectral approaches). 

Section \ref{sec:ModelIntro} introduces an illustrative model and demonstrates how it is defined on a multishape. In Section \ref{sec:2DChebClass} some of the important features of the 2DChebClass library are outlined, relevant for the MultiShape library. Section \ref{sec:ClassMethods} discusses the different class methods of the MultiShape implementation and Section \ref{sec:UserInterface} explains the user interface and demonstrates how to solve a problem of the form \eqref{eq:ProblemStatement1}, introduced in Section \ref{sec:ModelIntro}, on a multishape domain. 
\subsection{Model Equations} \label{sec:ModelIntro}
The MultiShape framework applies to a broad range of PDE models, some of which are introduced in Sections \ref{sec:Validation} and \ref{sec:ModelEquations}. In order to illustrate some of the complexities that are successfully addressed by the MultiShape framework, in this section we consider a system of non-linear, time-dependent (integro-)PDEs with (non-local) no-flux boundary conditions. This model is motivated by DDFT applications, which account for different species of interacting particles, and which were also used as illustrative examples in previous work, see \cite{GoddardPseudospectralCode1,OurPreprint2021}. For a system of $n_s$ species with densities $\rho_a$, $a = 1,...,n_s$, confined to a domain $\Omega \subset \mathbb{R}^2$, the general problem statement consists of a system of $n_s$ PDEs of the form
\begin{align}\label{eq:ProblemStatement1}
	\partial_t \rho_a(\vec x, t) &= -\nabla \cdot \mathfrak{ \vec j}_a(\{\rho_b\})    \qquad\qquad  \text{in } \Omega \times (0,T),\notag \\
	\mathfrak{ \vec j}_a(\{\rho_b\}) \cdot \vec n &= 0 \qquad\qquad \qquad  \qquad \hspace{0.2em} \ \ \ \text{on } \partial \Omega \times (0,T),\\
	\rho_a(\vec x, t)  &= \rho_{a,0}(\vec x) \ \qquad \qquad  \qquad \hspace{0.1em} \ \text{at } t = 0,\notag
\end{align}
where $\mathfrak{ \vec j}_a$ is a flux term and the notation $ \mathfrak{ \vec j}_a(\{\rho_b\})$, $b = 1,...,n_s$, denotes that the flux of species $a$ depends on the densities of all of the species. Note that the boundary condition is often referred to as `no-flux' and, in particular, enforces conservation of mass.  For DDFT applications, a common choice is to treat particle--particle interactions via the mean-field approximation \cite{DDFTReview}, resulting in a flux of the form 
\begin{align*}
	\mathfrak{ \vec j}_a(\{\rho_b\}) = -\nabla  \rho_a(\vec x, t)  - \rho_a(\vec x, t)  \nabla V_{\text{ext}}(\vec x, t)  - \sum_{b=1}^{n_s} \rho_a(\vec x, t) \int_\Omega \rho_b(\vec x',t) \nabla_{\vec x} V_{a,b}(|\vec x - \vec x'|) d\vec x ',
\end{align*}
where $V_{\text{ext}}$ is an external potential and $V_{a,b}$ are interaction potentials between species $a$ and $b$.

We solve the strong form of \eqref{eq:ProblemStatement1} on each of the elements. Matching conditions (patching) are applied at the intersection between two elements, and the PDE boundary conditions applied on the multishape domain boundary. The matching conditions between elements are
\begin{align*}
	\rho_{a,i}(\vec x, t)  &= \rho_{a,j}(\vec x, t) \qquad\qquad \quad \text{on } \partial \Omega_{i,j} \times (0,T), \notag\\
	-\mathfrak{ \vec j}_{a,i} \cdot \vec n_{i} &= \mathfrak{ \vec j}_{a,j} \cdot \vec n_{j}\ \ \ \ \  \quad \quad \qquad \text{on } \partial \Omega_{i,j} \times (0,T),\notag
\end{align*}
where $\partial \Omega_{i,j}$ denotes the intersection boundary between elements $e_{i}$ and $e_{j}$, $i,j = 1,...,n_e$, and quantities of the form $f_{a,i}$ denote the value of $f$ for species $a$ on element $i$.
In line with work on SEM for (classical) fluid dynamics problems, e.g., \cite[Chapter 5]{Canuto2007}, and since the system \eqref{eq:ProblemStatement1} is of interest in complex fluid dynamics, it is a natural choice to match the flux $\mathfrak{ \vec j}_a(\{\rho_b\})$ at the intersection of two elements, as compared to directly matching the first derivative of $\rho_a$. For problems of the form \eqref{eq:ProblemStatement1} these two approaches are essentially equivalent. The minus sign in the the flux matching condition arises due to the opposite orientation of (outward) normals for adjoining elements. Note that particle species are denoted by subscripts $a$, $b$, while elements have subscripts $i$, $j$.  

\subsection{2DChebClass} \label{sec:2DChebClass}
The code library 2DChebClass, developed in \cite{GoddardPseudospectralCode1}, supports the solution of systems of PDEs of the form \eqref{eq:ProblemStatement1}, as well as other PDE models, including the application of complicated boundary conditions, on individual elements. Note that in the 2DChebClass library these elements are referred to as shapes. 
The choice to solve the strong form of a PDE model and the use of a patching method in MultiShape is motivated by the fact that 2DChebClass provides a strong solution of a PDE on the individual elements. It is also advantageous for the end user, since the weak form of the PDE problem does not have to be derived in order to use this method; the PDE can be inputted directly as stated, e.g., in \eqref{eq:ProblemStatement1}. The choice of operator definitions and notation in MultiShape ensures that writing a PDE and the corresponding boundary conditions
can be done in a `natural' way within the framework.  Additionally, as will be demonstrated in Section \ref{sec:NumericalMethods}, extra care has been taken to ensure that the implementation of (complicated, non-linear, non-local) boundary conditions, as well as matching conditions between the elements, is straightforward for the user and comparable to the 2DChebClass library. We briefly outline some of the features of 2DChebClass that are carried forward to the MultiShape library, and which are explained in detail in \cite{GoddardPseudospectralCode1}. In particular, the MultiShape implementation aims to be intuitive for both 2DChebClass and new users.

Whilst pseudospectral methods are well suited to periodic problems, such domains are not directly applicable in a spectral element setting.
Hence, here we restrict ourselves here to discussing the case in which each dimension of a square, $[-1,1]$, is discretized using Gauss--Lobatto--Chebyshev collocation points, defined as
\begin{align*}
	x_i = \cos \left(\frac{i \pi}{2}\right),
\end{align*}
for $i = 0,1,...,N-1$. In the two-dimensional domain $[-1,1]^2$, Kronecker products are used to create a two-dimensional grid in the standard way \cite[Chapter 7]{bibTrefethen}. The domain $[-1,1]^2$ is referred to as the computational domain and, using a bi-linear map, it can be associated with a quadrilateral (in the case of 2DChebClass, only a rectangular box is provided by default) or wedge element, which is referred to as the physical domain. Note that for wedge elements, both the computational and physical domain are defined in polar coordinates.
The interpolation polynomial on an element is defined using barycentric Lagrange interpolation \cite{bibTrefethenBerrut1}. In the 2DChebClass implementation first and second (and, optionally, higher) order differentiation matrices are constructed automatically on the discretized shape, as well as gradient, divergence, and Laplacian operators for two-dimensional problems. These operators are constructed on the computational domain $[-1,1]^2$ and then mapped onto the physical domain using the aforementioned (inverse of the) bi-linear map, making it easy to define operators on different physical domains \cite[Chapter 22]{Boyd1}. For the evaluation of integrals, Clenshaw--Curtis quadrature \cite{ClenCurt1} is applied. For further details on the pseudospectral implementation in 2DChebClass, we refer to \cite{GoddardPseudospectralCode1}.

\subsection{Class Methods} \label{sec:ClassMethods}

In this section, technical aspects of the MultiShape implementation, as a class in {\scshape Matlab}, are discussed, such as the construction of differentiation and interpolation matrices, integration vectors, convolution operators, and the application of intersection conditions, all of which are defined as class methods. While many quantities can be computed on the individual elements first, a key part of the implementation is that all operators can be applied to the whole multishape domain simultaneously, which is described in the following. Note that, since quadrilateral and wedge elements are defined in two different coordinate systems, Cartesian and polar, the corresponding operators are defined in terms of two different coordinate systems as well.

In line with the philosophy of the MultiShape code, it draws on the 2DChebClass implementation for single shapes. Discretizing a function $f$ on a multishape domain is done by discretizing $f$ on each element individually through 2DChebClass, and stacking the resulting vectors into a vector of size $M \times 1$. Here $M = \sum_{i=1}^{n_e} N^{(i)}_1 N^{(i)}_2$, where $N^{(i)}_1$, $N^{(i)}_2$ are the number of points for element $e_i$ in the the $x_1$ and $x_2$ directions, respectively, and $n_e$ the number of elements in the multishape. The resulting discretized function is
\begin{align*}
	\mathbf f_{M \times 1} = \begin{pmatrix}
		\mathbf f_{1} \\ \mathbf f_{2} \\ \vdots \\ \mathbf f_{n_e}
	\end{pmatrix},
\end{align*} 
where $\mathbf f_{i}$ is a $N_1^{(i)}N_2^{(i)} \times 1$ vector, for $i = 1,..., n_e$.
The same approach is taken when constructing the integration vector for the multishape, which is a row vector of size $1 \times M$. 

The gradient, divergence, and Laplacian operators for a multishape are constructed by computing the operators on each individual element, using 2DChebClass, and arranging them in a block diagonal matrix, as demonstrated below for an operator $D$ acting on a function (vector) $\mathbf f$
\begin{align*}
	D \mathbf f = \begin{pmatrix}
		D_{1} & & & &\\
		& D_{2} & &\\
		& & \ddots &\\
		& & & D_{n_e}
	\end{pmatrix}
	\begin{pmatrix}
		\mathbf f_{1} \\ \mathbf f_{2} \\ \vdots \\ \mathbf f_{n_e}
	\end{pmatrix}.
\end{align*}
Note that empty blocks in the above denote zeros.
Computations can then be carried out directly on the multishape, instead of on each element individually, since the stacking ensures that the order of elements is preserved. A matrix $D_{i}$ corresponding to element $e_{i}$ is applied only to $\mathbf f_{i}$, corresponding to the same element, within the matrix--vector multiplication. Note that $D_{i}$ and $D_{j}$ could be defined in two different coordinate systems, depending on whether, for example, $e_{i}$ and $e_{j}$ are quadrilateral or wedge elements.

Interpolating a function $f$ from one grid onto another is another crucial function of the MultiShape implementation. Accurate interpolation is relevant for plotting of results, computation of convolutions with finite support, and evaluation of quantities on subdomains, as well as comparing results evaluated on different grids. Many of these applications only require the interpolation matrices in computational space $[-1,1]^2$ for each element individually, which can simply be stacked, as illustrated above for the matrix $D$, and applied to a vector $\mathbf f$. 
However, mapping a function $f$ from a set of points $\mathbf x_{\rm in}$ onto another arbitrary set of points $\mathbf x_{\rm out}$, which is defined in physical (not computational) space, is more complicated than stacking the individual interpolation matrices and applying these to $ \mathbf f(\mathbf x_{\rm in})$. The reason for this is that it is not clear a-priori in which element of the multishape the different points within $\mathbf x_{\rm out}$ lie. The code is equipped to determine this automatically for the user, by identifying which of the $\mathbf x_{\rm out}$ lie in which shape. Once this is established, the interpolation matrix for each element $e_{i}$ can be used to interpolate onto the set of points in that element, $\mathbf x_{{i},\rm{out}}$.

In order to apply boundary and matching conditions we need to identify the boundaries of the multishape and the intersections between individual elements. To aid illustration, Figure \ref{fig:MS1} displays a multishape, consisting of two elements, with marked intersection boundary, multishape boundary, and normal vectors to the boundary. 2DChebClass is able to automatically compute the boundaries of a single element $e_{i}$, $i = 1,..., n_e$. It furthermore distinguishes the four parts, or faces, of the boundary, i.e., `top', `bottom', `left', and `right' for a quadrilateral; inner and outer radial and maximum and minimum angular boundaries for a wedge. We denote these faces by $f_{i,k}$, $k = 1,2,3,4$. In order to find the intersections between two elements, we have to check whether the points on face $f_{i,k}$ of element $e_{i}$ are equal to those of face $f_{j,l}$ of element $e_{j}$ for each possible combination of elements and faces. The ordering of the points lying on face $f_{i,k}$ can also be the reverse of the order of points on $f_{j,l}$, which MultiShape checks automatically. 
Then, the multishape boundary can be found by subtracting those intersection boundaries from the union of boundaries of the individual elements, leaving only those parts of the element boundaries that are not adjacent to another element. In this way, the patching method is implemented -- note that, in principle, it is possible to modify this approach and admit the intersecting faces of two elements to have different numbers of points, then use MultiShape's interpolation function to match these. However, since a similar number of points is generally expected in neighbouring elements, this additional complexity is omitted in the current implementation.

\begin{figure}[h]
	\centering
	\includegraphics[scale=0.17]{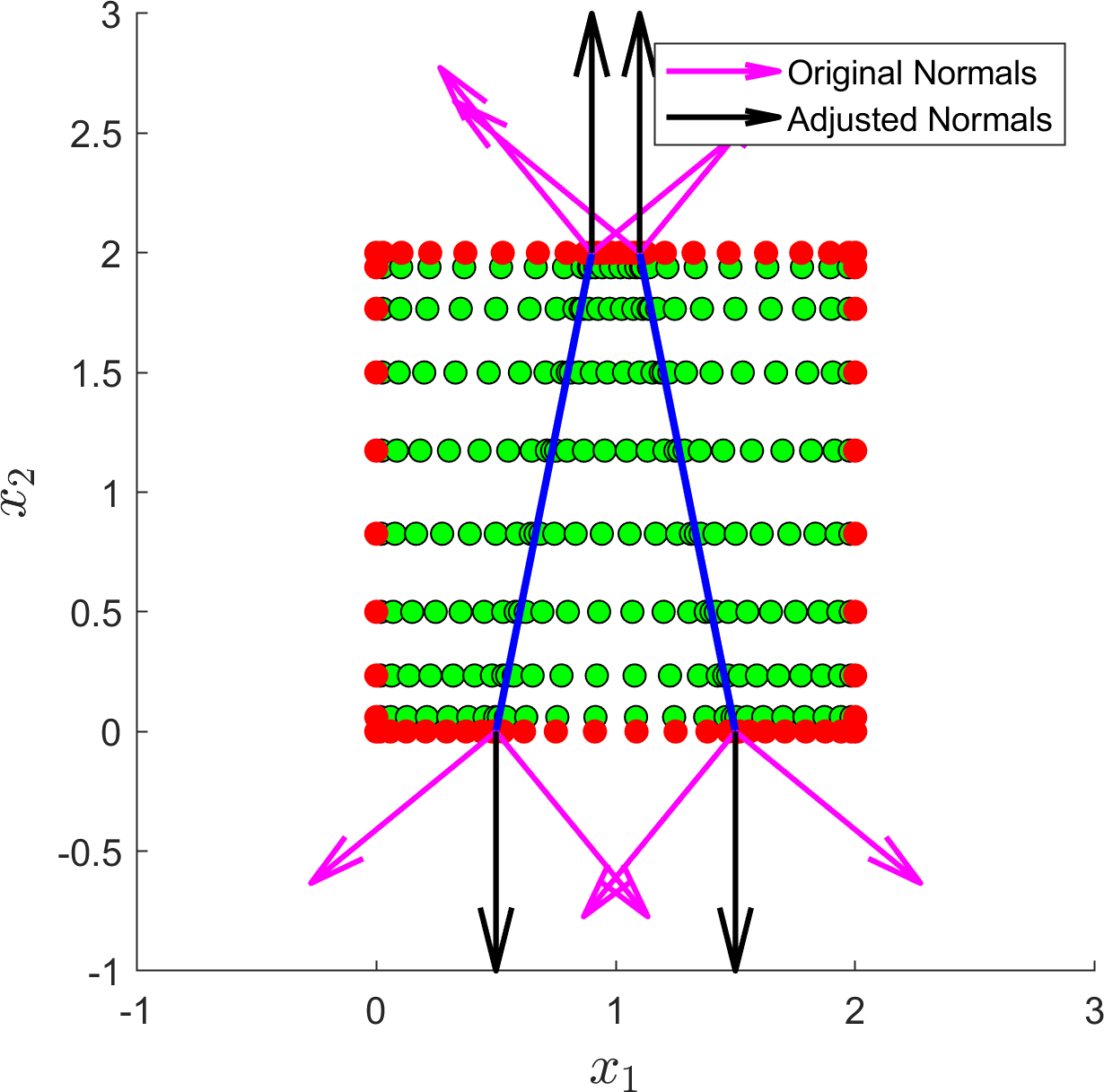}
	\caption{Overriding the normal vectors at intersection corners, which are automatically constructed by MultiShape, to ensure they are defined in the normal direction to the multishape face on which the intersection corner is located. } 
	\label{fig:NormalOverride}
\end{figure}

Constructing the normal vectors on the boundary of the multishape, which are needed to compute no-flux boundary conditions and intra-element matching conditions, is a little more complicated. The starting point is again to consider the normals for each individual element, provided by 2DChebClass, and in particular those lying on the boundary of the multishape. However, at the intersection corners between two elements, these normal vectors are not defined consistently. In this case the code considers the normal vectors on the multishape boundary to both sides of the intersection corner and averages them. In order to do this correctly, all normals first have to be converted into Cartesian coordinates, averaged, and converted back into polar coordinates where appropriate. When the discretization of the multishape is more complicated, taking the average of adjacent normals may not be sufficient, and the resulting outward normal may not be sensible. Since this can be very specific to the discretization at hand, the option to override any of the normal vectors manually is given. Such a case is demonstrated in Figure \ref{fig:NormalOverride}.

The computation of convolution integrals is the final part of the MultiShape implementation that we discuss here. Since computing a convolution integral is a global operation, we cannot rely on reusing the computations on individual elements. However, the general steps are equivalent to the approach for a single element and it is therefore recommended to refer to \cite{GoddardPseudospectralCode1} for details. A brief outline of constructing the convolution matrix is given below. A convolution integral is defined as
\begin{align*}
	\left(\rho \star \chi \right) (\vec y) = \int_\Omega \chi \left( \vec y - \vec{ z}\right) \rho \left(\vec{z}\right) d \vec{z},
\end{align*}
for a density $\rho$ and a weight function, or kernel, $\chi$. In the system \eqref{eq:ProblemStatement1}, $\chi = \nabla V_{a,b}$. Given the integration vector, Int, for the multishape and defining $ \vec y - \vec z: = \vec x$, the convolution matrix can be computed by
\begin{align*}
	\text{Conv}_{m,n} = \text{Int}_n \ \ \chi\left({x}_1^{(m)} - {x}_1^{(n)},{x}_2^{(m)} - {x}_2^{(n)}\right),
\end{align*}
where the subscripts $1$ and $2$ refer to the two dimensions of the domain, and $m, n \in \{1,...,M\}$, with $M$ the total number of points in the multishape.
If the function $\chi$ only takes one input, we consider 
\begin{align*}
	\text{Conv}_{m,n} = \text{Int}_n \ \ \chi\left(\|d_1^{(m,n)}, d_2^{(m,n)}\|_2\right),
\end{align*}
where $d_1^{(m,n)} := {x}_1^{(m)} - {x}_1^{(n)}$, $d_2^{(m,n)} := {x}_2^{(m)} - {x}_2^{(n)}$, and $\| \cdot \|_2$ denotes the standard Euclidean distance.
The convolution integral is now defined as a matrix--vector multiplication, by applying the matrix Conv to a density vector $\rho$. 
The convolution matrix, which is independent of the density $\rho$, can be applied to different functions, without needing to be recomputed. In particular, significant computation time is saved in outer routines, such as ODE solvers and optimization algorithms, in which the same convolution matrix can be applied to an updated $\rho$ at each step of the sub-routine. 
This is contrasted by the fast Fourier transform (FFT) approach for solving convolution integrals \cite{roth2010fundamental}. Here, at each step in time a Fourier transform of $\rho$ is taken, and the convolution theorem is applied to yield
\begin{align*}
	\left(\rho \star \chi \right) (\vec y) =  F^{-1} \{  F \{\rho\} \cdot F \{\chi\}\}(\vec y).
\end{align*}
It is clear that the present approach does not exhibit these complications, and therefore can lead to noticeable computational advantages when used within a sub-routine. 

\subsection{User Interface} \label{sec:UserInterface}
This section illustrates how a multishape is set up by the user and describes the information that needs to be provided in order to solve PDE problem \eqref{eq:ProblemStatement1} on the multishape using the DAE solver \texttt{ode15s} in {\scshape Matlab}.  

In order to set up a multishape, each element has to be defined, as well as an order in which they are arranged, from $1$ to $n_e$, although this ordering
is essentially arbitrary. For each element the number of discretization points, $N_1$ in the $x_1$-direction, $N_2$ in the $x_2$-direction, needs to be specified. To uniquely determine a shape, for a quadrilateral element it suffices to specify the corners; a wedge element is defined in polar coordinates, by specifying an inner and outer radius, a minimum and maximum angle, and the location of the origin. 

In PDE problems such as \eqref{eq:ProblemStatement1}, matching conditions (see Section \ref{sec:ModelIntro}) are applied at the intersections between two elements $e_{i}$ and $e_{j}$. However other conditions can be applied at the intersection boundary between two elements. For example, the MultiShape library additionally supports the application of no-flux boundary conditions between elements, modelling a hard wall, and the inclusion of other matching options into the code library is straightforward. Which condition to apply on each intersection between elements $e_{i}$ and $e_{j}$, can be specified simply by setting a flag. The code demonstrating how to set up a multishape, and provide the necessary geometrical and numerical information, is found in Listing \ref{code:MSSetup}. The resulting multishape can be seen in Figure \ref{fig:MS1}. Note that, in order to create a functioning multishape, the number of points of the neighbouring faces of two elements have to match, as discussed in Section \ref{sec:ClassMethods}. Once the multishape is set up, required operators and properties, such as differentiation matrices or the location of boundary points, can be simply extracted from the multishape structure, as demonstrated in Listing \ref{code:OpSetUp}. 

Using MultiShape to discretize a PDE and apply boundary and intersection conditions is straightforward. A discretized version of problem \eqref{eq:ProblemStatement1} can be defined by stacking $n_s$ discretized equations of the form
\begin{align*}
	M \boldsymbol{\rho}_{a}'(t) = \mathbf{f}(\{\boldsymbol{\rho}_{b}(t) \},t),
\end{align*}
where the mass matrix $M$ for a collocation method is diagonal, and $\mathbf{f}$ is the discretization of the right hand side of the PDE. The boundary and matching conditions are applied by setting the relevant (diagonal) entries in the mass matrix to zero. Discretizing the right-hand side of the PDE, $\mathbf{f}(\{\boldsymbol{\rho}_{b}(t)\},t)$, is illustrated in Listing \ref{code:PDEdisc}. The differentiation and convolution matrices are constructed automatically by the MultiShape library and can be applied to a vector $\boldsymbol{\rho}_{a}$ using matrix--vector multiplication. The intersection conditions are applied via an inbuilt function, which uses the intersection boundaries that are automatically identified during construction of the multishape. The quantities to be provided are the solution, $\boldsymbol \rho_a$, and the flux, $\mathbf j (\{\boldsymbol \rho_b\})$, that are to be matched. Similarly, the code automatically identifies the boundary of the multishape, so that the boundary conditions of the PDE can be applied in the same way as for single shapes in 2DChebClass. Note that the non-local no-flux boundary conditions of problem \eqref{eq:ProblemStatement1} are applied in Listing \ref{code:PDEdisc} in a single line of code (line 23), demonstrating the ease of implementation of complicated boundary conditions. Additionally, this line can be replaced to apply other boundary conditions, such as Dirichlet or Robin conditions. 
Note that several `helper functions' are defined within MultiShape, which support straightforward implementation. These include masks to access individual elements of the multishape, computation of inner products between two vectors, and the application of intersection boundary conditions, as discussed in Section \ref{sec:ClassMethods}. 

In order to solve the model problem \eqref{eq:ProblemStatement1} using the {\scshape Matlab} DAE solver \texttt{ode15s}, one needs to specify the right-hand-side of the discretized PDE as illustrated above, an initial condition, a time interval, and a mass matrix, as well as the solver tolerances. This is illustrated in Listing \ref{code:PDEsol}. Note that the codes in Listings \ref{code:OpSetUp}, \ref{code:PDEdisc}, and \ref{code:PDEsol} are independent of the specific multishape defined in Listing \ref{code:MSSetup}. If one were to compute the same PDE system on a different multishape, the only part of the code that needed to be amended would be in Listing \ref{code:MSSetup}, demonstrating the modularity of the code.

\begin{lstlisting}[frame=single, label=code:MSSetup,caption=Setting up a multishape]  
	
% Number of points (note N2x1 = N1x1 since the numbers of points need 
% to match at intersection boundaries)
N1x1 = 12; N1x2 = 15; N2x1 = N1x1; N2x2 = 15;

% For a quadrilateral:
% - specify corners and number of points
x1Min = 0; x2Min = 0; x1Max = 3; x2Max = 3;
geom1.Y = [x1Min,x2Min;x1Min,x2Max;x1Max,x2Max;x1Max,x2Min];
geom1.N = [N1x1;N1x2];
% - assign the geometry to be shape 1; specify that it is a quadrilateral
shapes(1).geom = geom1; shapes(1).shape = 'Quadrilateral';

% For a wedge:
% - specify inner/outer radius, minimum/maximum angle, and origin 
rMin = 1; rMax = 4; thetaMin = 0; thetaMax = pi;
geom2.R_in = rMin; geom2.R_out = rMax;
geom2.th1 = thetaMin; geom2.th2 = thetaMax;
geom2.Origin = [4,3];
geom2.N = [N2x1;N2x2];
% - assign the geometry to be shape 2; specify that it is a wedge
shapes(2).geom = geom2; shapes(2).shape = 'Wedge';

% Match intersections between shape 1 and shape 2
BCs = struct; BCs(1,2).function = 'BCmatch';
% Construct multishape, given the shapes
MS = MultiShape(shapes);
\end{lstlisting}

\begin{lstlisting}[label=code:OpSetUp, caption=Extracting operators]
% Extract pre-computed operators from multishape MS: this improves
% performance of the ODE solve, since accessing structures in MATLAB
% is comparatively slow
grad = MS.Diff.grad; div = MS.Diff.div; Lap = MS.Diff.Lap; Int = MS.Int;
bound = MS.Ind.bound; normal = MS.Ind.normal;
intersections = MS.Ind.intersections;
x1 = MS.Pts.y1_kv; x2 = MS.Pts.y2_kv;

% Compute convolution matrix, for a kernel defined by a function V2Fun
Conv = MS.ComputeConvolutionMatrix(@V2Fun,true); 
\end{lstlisting}

\begin{lstlisting}[label=code:PDEdisc, caption=Discretizing the right-hand side of the PDE system]
% Right-hand side function, including boundary & intersection conditions
function drhodt = drho_dt(t,rhoab)
	% Extract the individual species
	rhoa = rhoab(1:end/2,:); rhob = rhoab(end/2+1:end,:);
	
	% Define discretized PDE using known external potential vector
	% (Vext), differentiation matrices (grad, div), and convolution
	% matrices between species a and b (Convaa, Convab, Convbb, Convba)
	rhoaflux = Flux(rhoa,rhob,Convaa,Convab);       
	drhoadt  = -div*rhoaflux;
	rhobflux = Flux(rhob,rhoa,Convbb,Convba);        
	drhobdt  = -div*rhobflux;
	
	% Apply matching conditions at intersections
	dataa.rho = rhoa; dataa.flux = rhoaflux;
	drhoadt = MS1.ApplyIntersectionBCs(drhoadt,dataa,BCs);
	datab.rho = rhob; datab.flux = rhobflux;
	drhobdt = MS1.ApplyIntersectionBCs(drhobdt,datab,BCs);
	
	% Apply no flux boundary conditions at boundary
	drhoadt(bound) = normal*rhoaflux; drhobdt(bound) = normal*rhobflux;
	drhodt = [drhoadt;drhobdt];
end
function rhoFlux = Flux(rhoa,rhob,Convaa,Convab)
	% Stack rho - this is now a size 2M x 1 vector
	rhoa2 = MS.MakeVector(rhoa,rhoa); rhob2 = MS.MakeVector(rhob,rhob);
	
	rhoFlux = -(grad*rhoa + rhoa2.*(grad*Vext) + rhoa2.*(grad*(Convaa*rhoa)) + rhoa2.*(grad*(Convab*rhob)));
end
	
\end{lstlisting}

\begin{lstlisting}[label=code:PDEsol, caption=Solving a PDE system on a multishape using {\scshape Matlab}'s DAE solver]
% Solving a system of PDEs using MATLAB's inbuilt DAE solver
% Set the mass matrix 
mM = ones(size(x1));   % vector of ones, of size M x 1
                       % M = length of discretized spatial domain
mM(bound) = 0;         % set to zero at boundaries
mM(intersections) = 0; % set to zero at intersections
MIn = diag([mM;mM]);   % diagonal matrix to account for system of 2 PDEs

% Set options for ODE solver: relative & absolute tolerance, mass matrix
opts = odeset('RelTol',10^-9,'AbsTol',10^-9,'Mass',MIn);
% Solve ODE, given the following:
% - vector of time points comp_times
% - initial condition rhoab_ic 
% - right-hand side function 'drho_dt', computed in Listing 3
[outTimes, rhoab_t] = ode15s(@drho_dt,comp_times,rhoab_ic,opts);  

\end{lstlisting}

\section{Validation of MultiShape Methodology} \label{sec:Validation}
This section contains validation tests for our MultiShape implementation, investigating the accuracy of the differentiation, integration, interpolation, and convolution operations. We further validate our methodology by computing the solution to a Poisson equation with Dirichlet boundary conditions. All of these tests have been carried out by comparing the numerical results to analytic solutions. In all tests exponential convergence is observed, until a certain threshold is reached, which is mostly dictated by machine precision, and to a lesser extent by compounding precision or interpolation errors.

Most errors in this section are calculated using a relative numerical $L^2$ norm defined as
\begin{align}\label{eq:ErrorMeasureL2}
	\mathcal E = \frac{||f_{\text{Num}} - f_{\text{Ex}}||_{L_2}}{||f_{\text{Ex}}||_{L_2} + 10^{-10}},
\end{align}
where $f_{\text{Num}} $ and $f_{\text{Ex}}$ correspond to numerical and exact solutions, and the additional (arbitrary) small term in the denominator is added to avoid the possibility of division by zero -- the results do not change significantly for other reasonable values. The interpolation error is calculated in the same way but replacing the $L^2$ norm by a (vector) $\ell_\infty$ norm. This is because the error is calculated on a uniform grid and not a Chebyshev grid, so that computation of the $L^2$ norm would require numerical integration in uniform space, which would introduce larger integration errors than the Clenshaw--Curtis integration formula \cite{ClenCurt1} implemented in MultiShape. Since integration results in a single number, for such computations, the $L^2$ norm in \eqref{eq:ErrorMeasureL2} is replaced by taking absolute values, with $f_{\text{Num}}$ and $f_{\text{Ex}}$ the numerical and exact integrals, respectively. 

\subsection{Validation of Discretized Operators}\label{sec:ValidationOP}
\begin{figure}[h]
	\centering
	\includegraphics[scale=0.1]{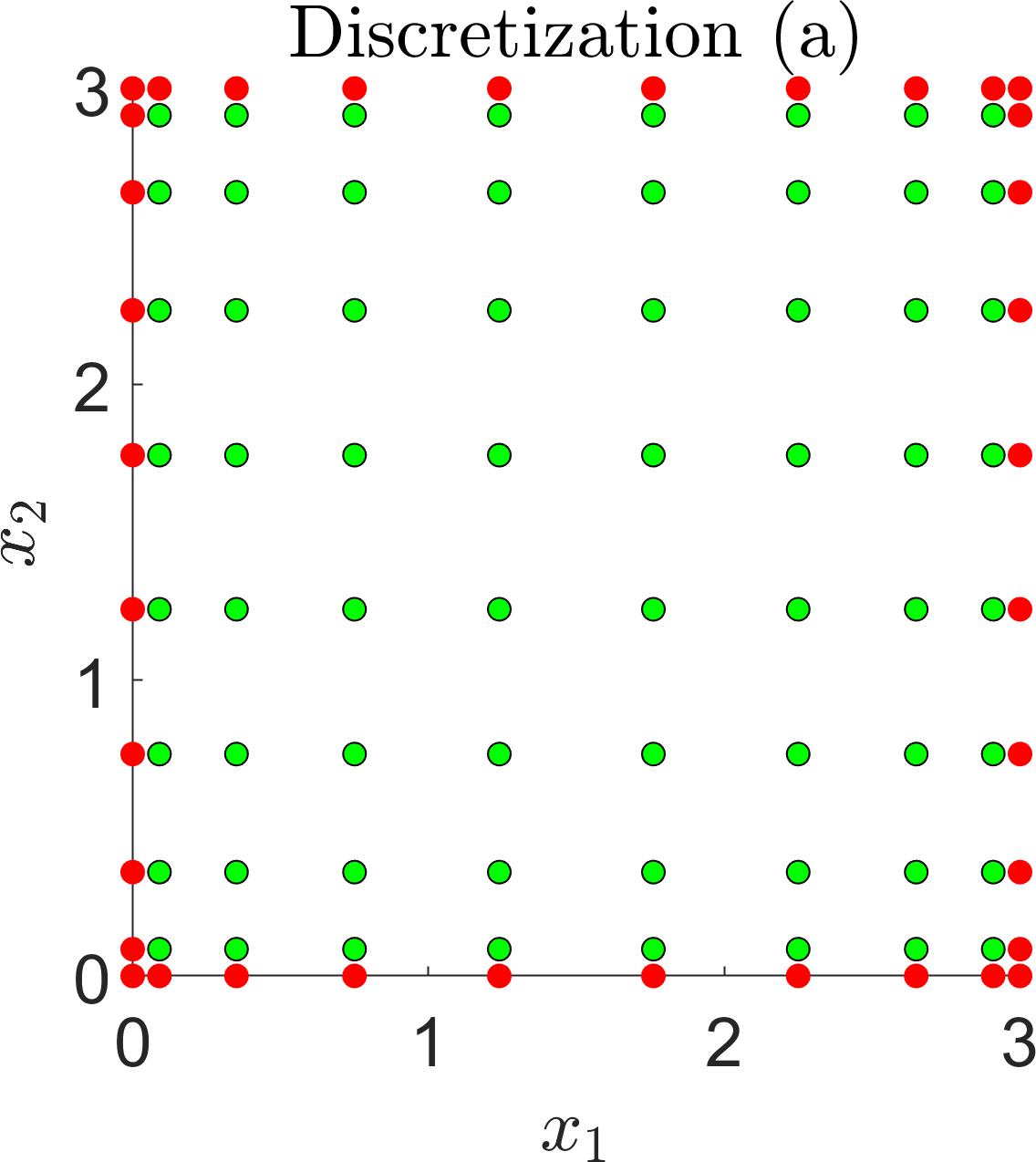}
	\includegraphics[scale=0.1]{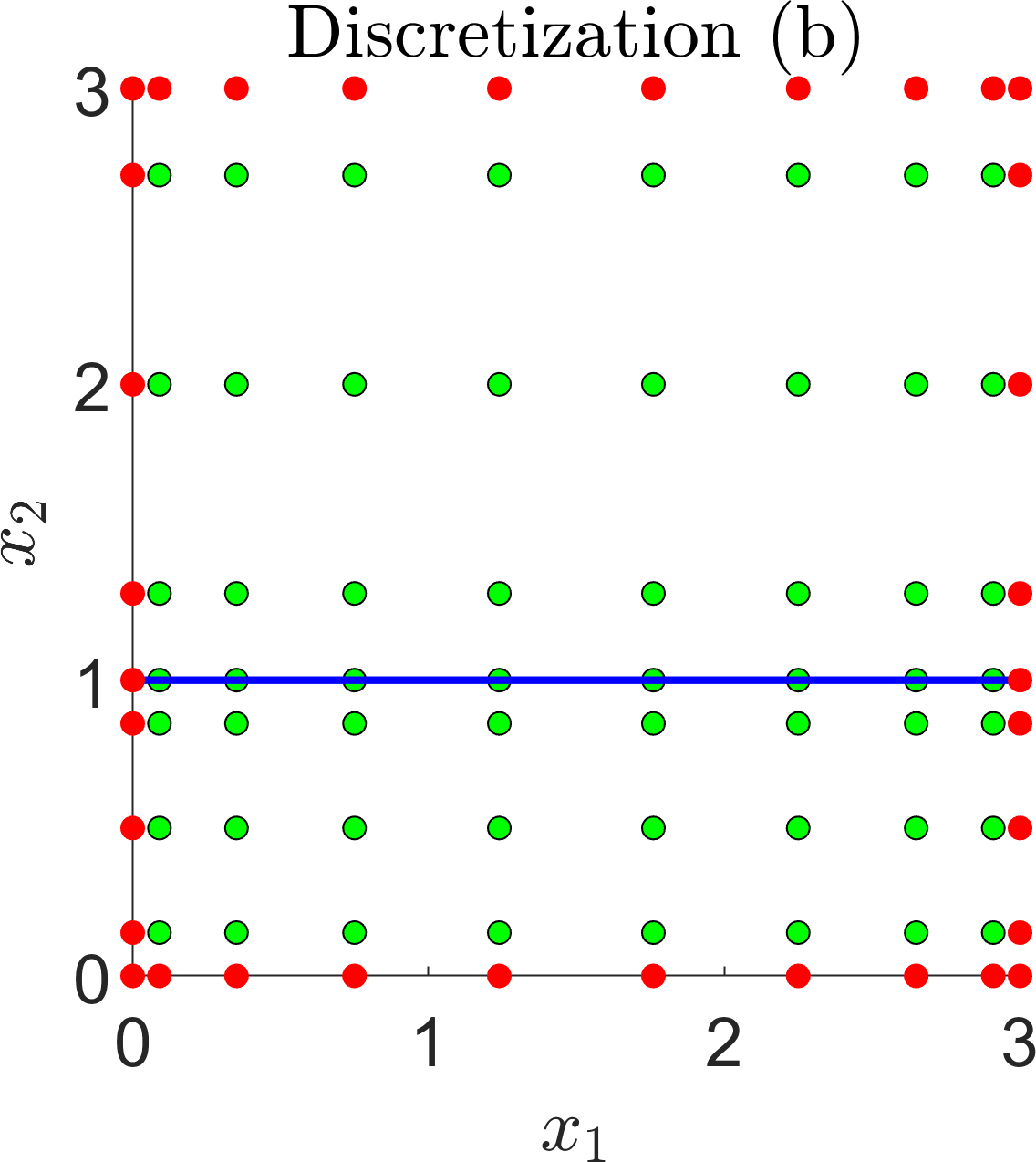}\\
	\includegraphics[scale=0.1]{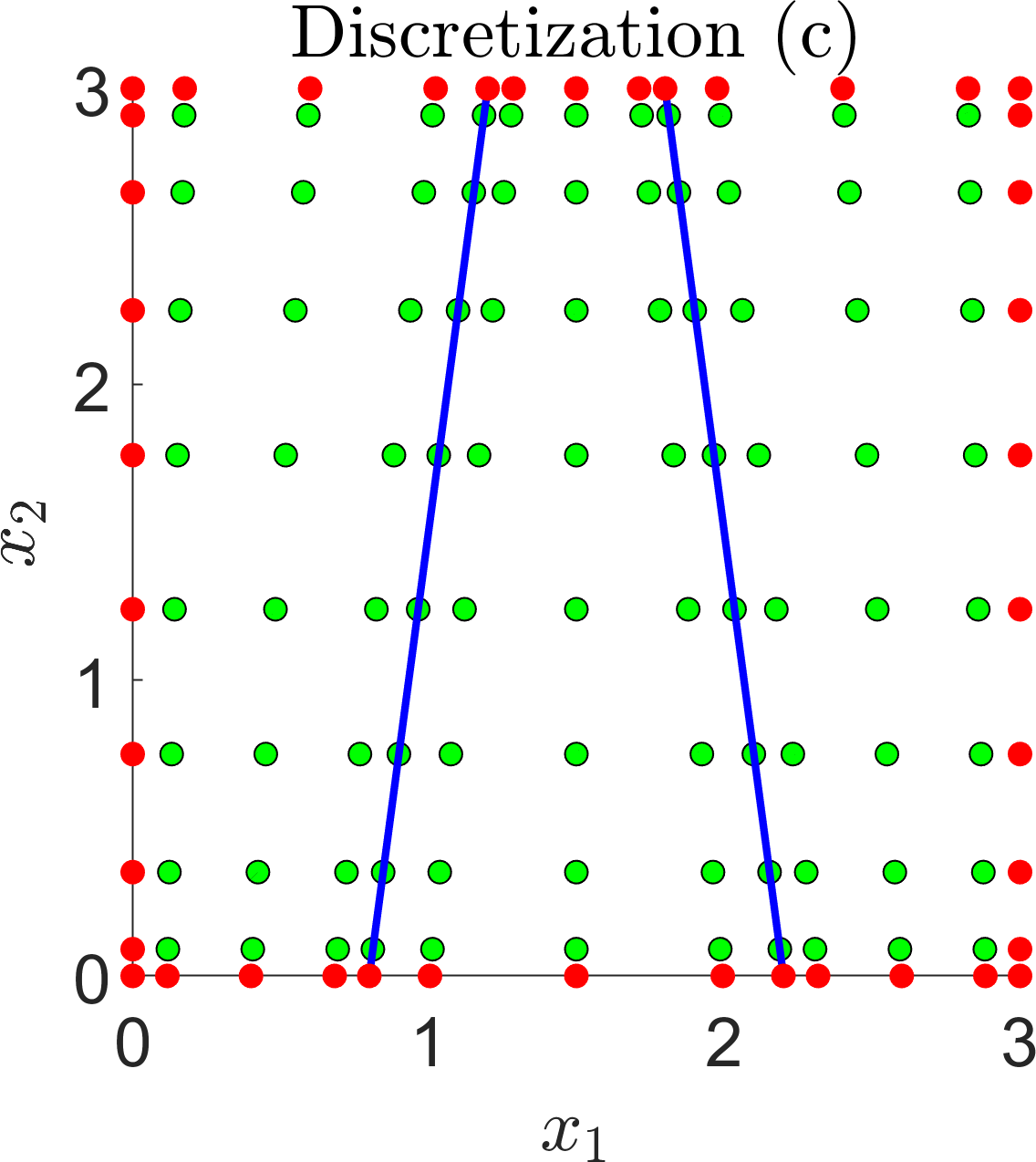}
	\includegraphics[scale=0.1]{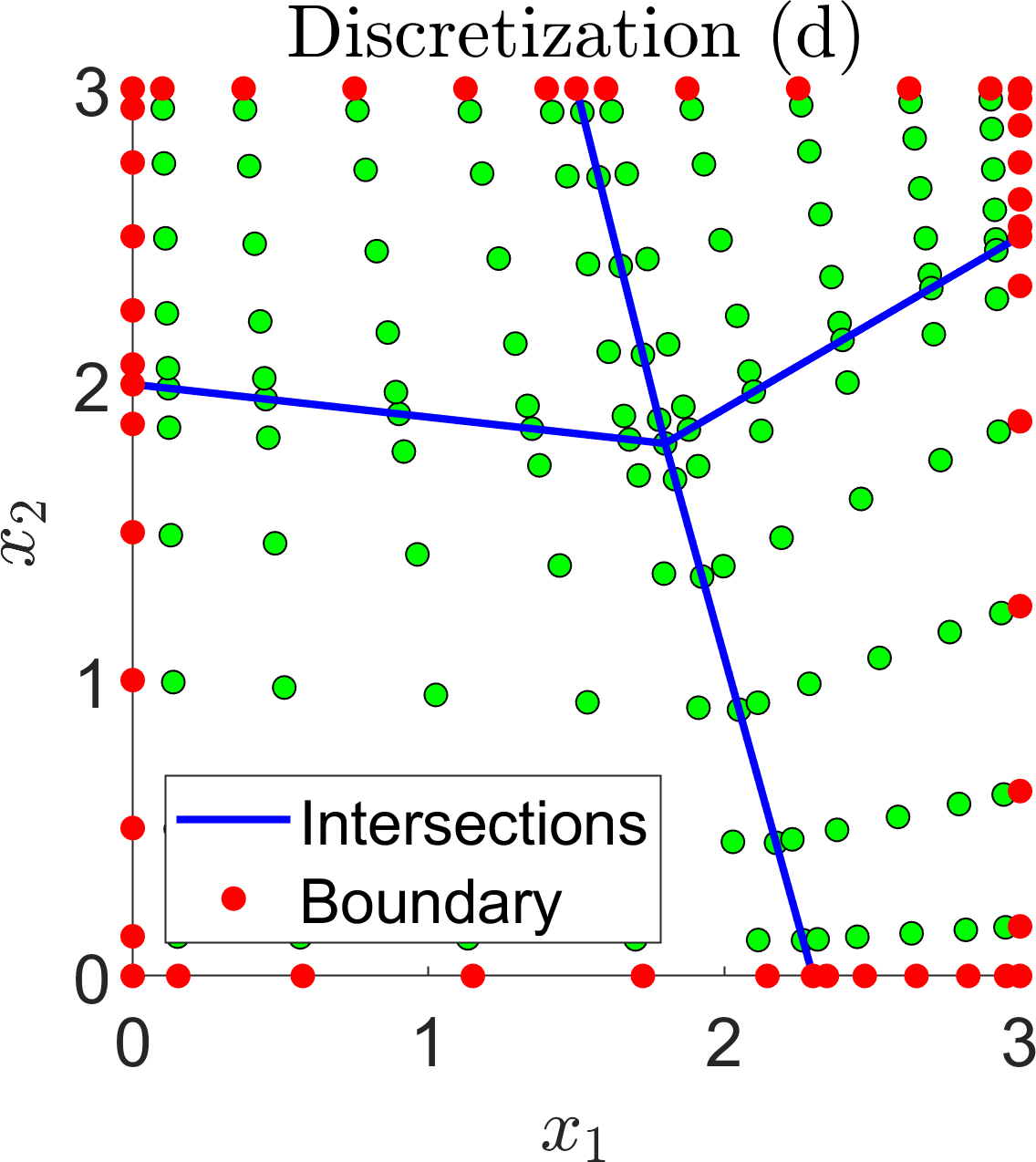}
	\caption{Different multishape discretizations of a box (a--d).}
	\label{fig:BoxSections}
\end{figure}

\begin{figure}[h]
	\centering
	\includegraphics[scale=0.1]{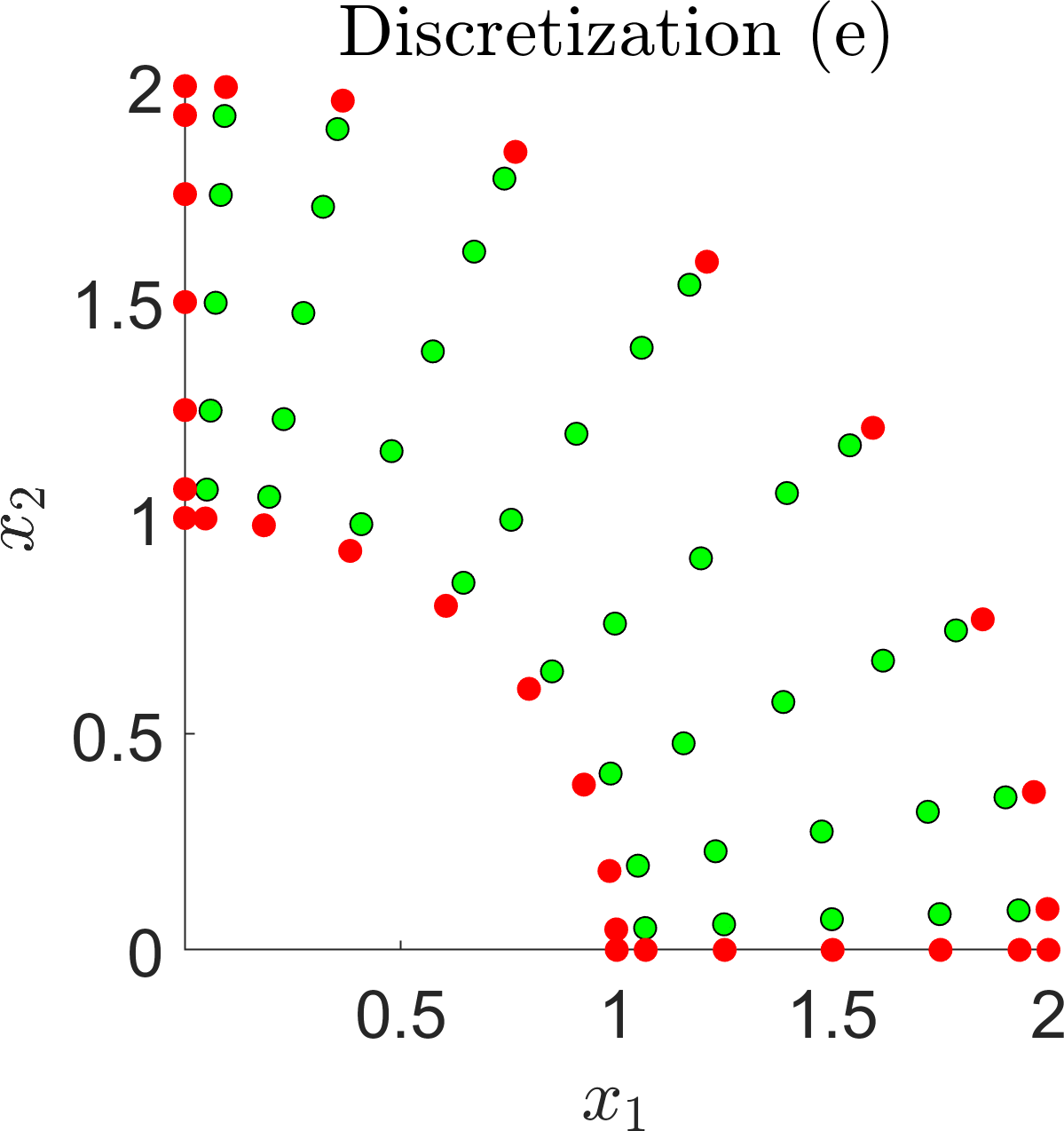}
	\includegraphics[scale=0.1]{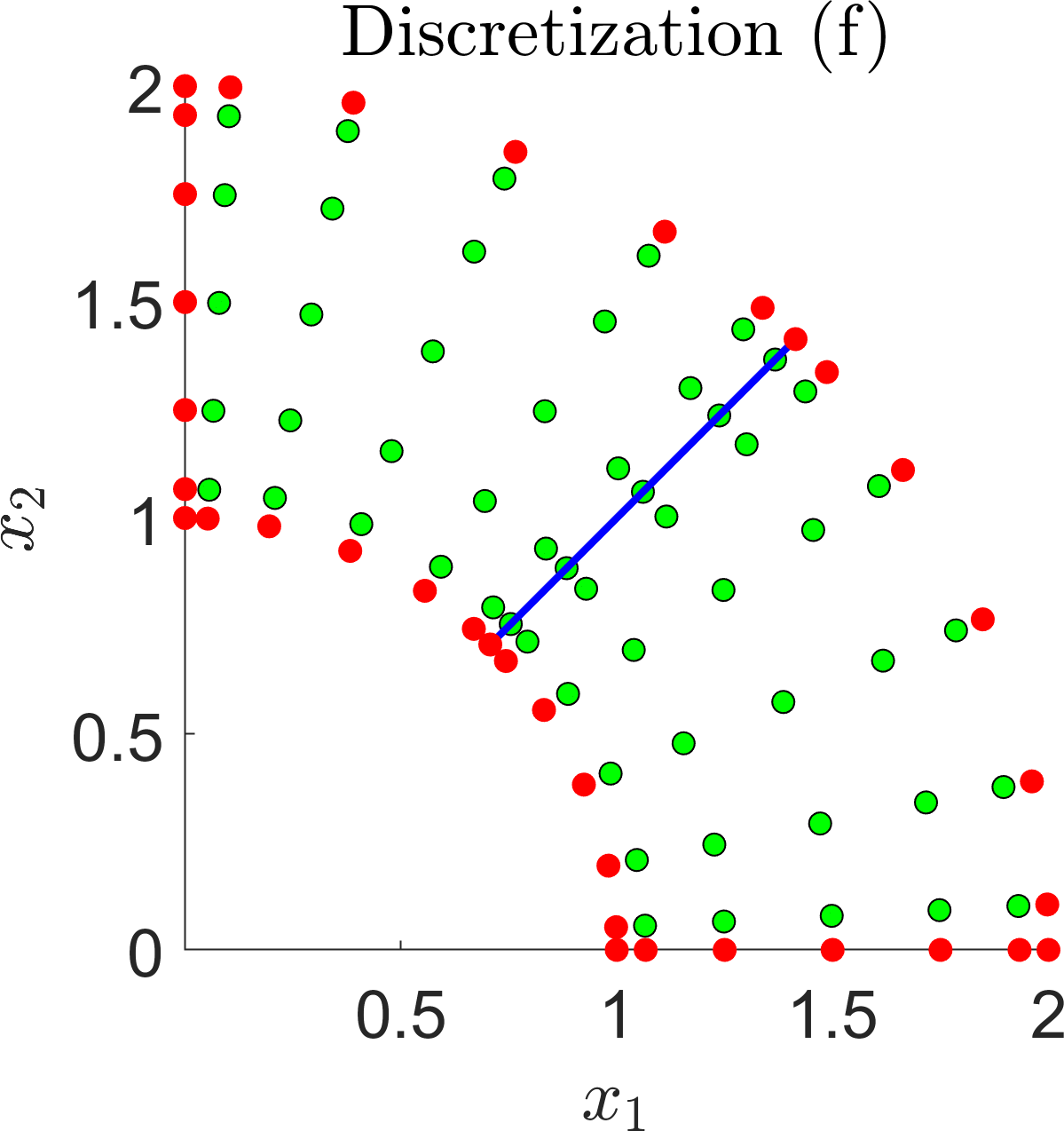}\\
	\includegraphics[scale=0.1]{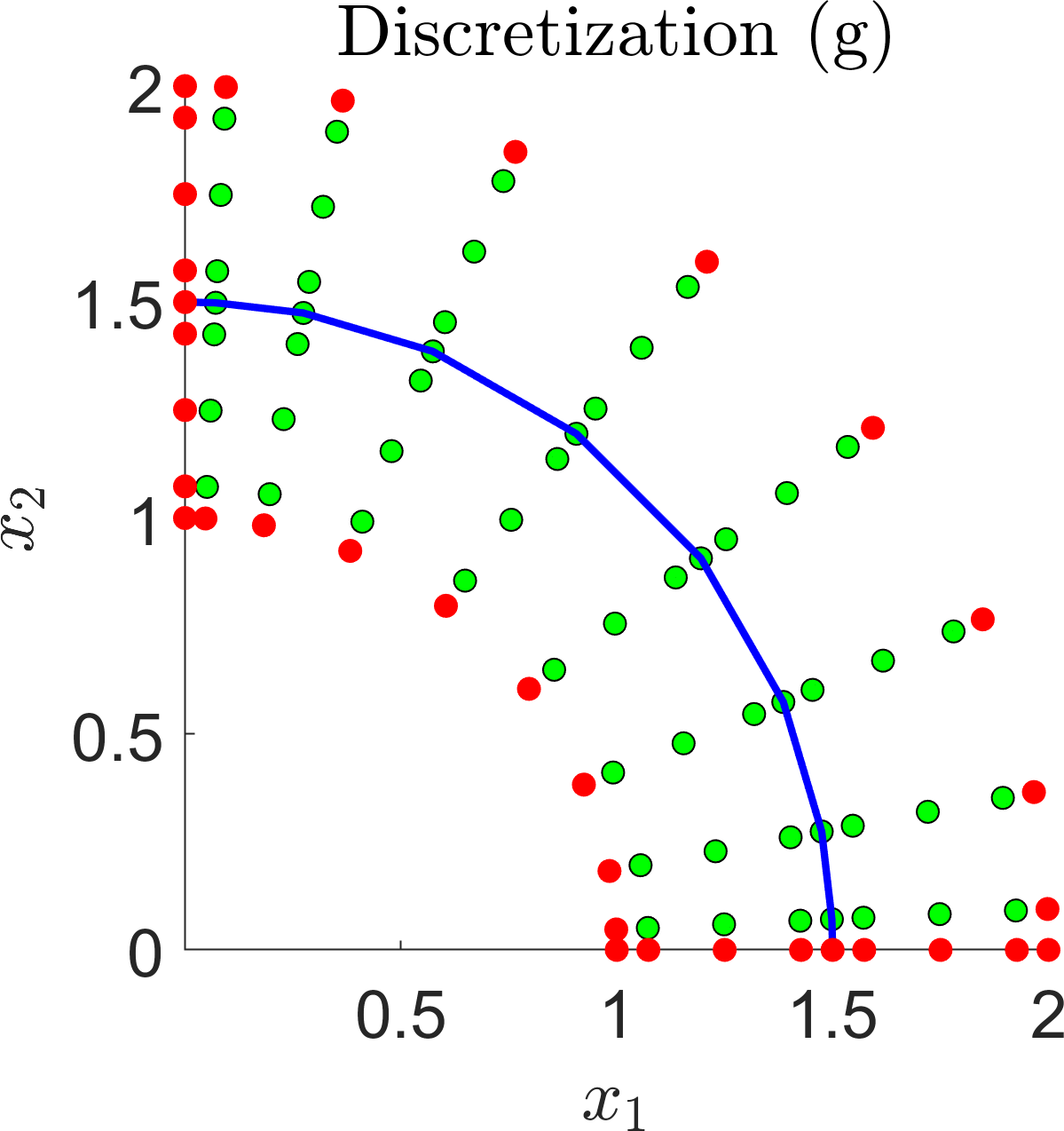}
	\includegraphics[scale=0.1]{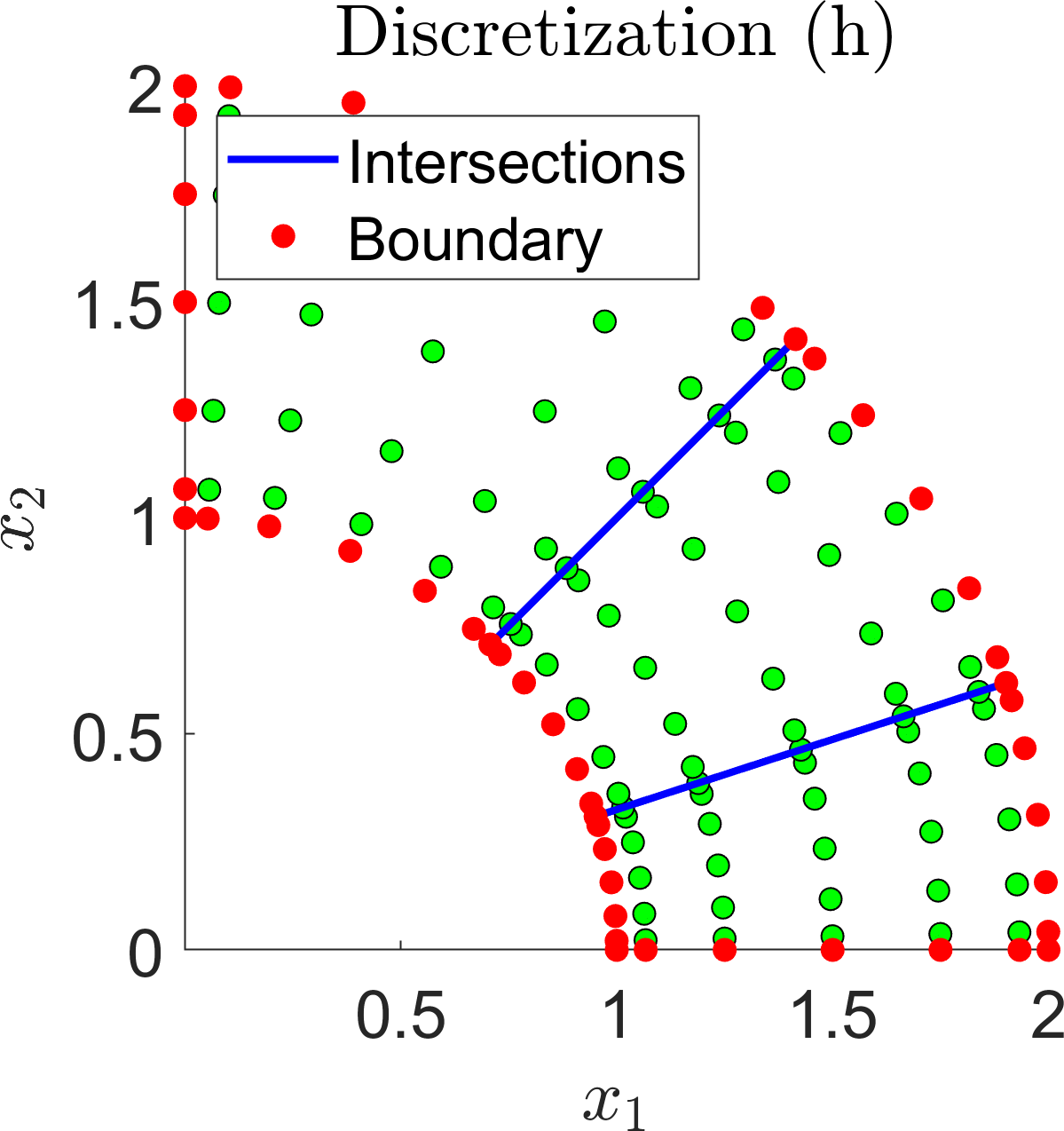}
	\caption{Different multishape discretizations of a wedge (e--h).} 
	\label{fig:WedgeSections}
\end{figure}

The multishapes considered in this validation section are discretizations of boxes and wedges. This is so that comparisons can be made with the 2DChebClass implementation  for a single shape. The different discretizations considered for a variety of test cases can be seen in Figures \ref{fig:BoxSections} and \ref{fig:WedgeSections}. The four discretizations of a box in Figure \ref{fig:BoxSections} include a reference multishape, discretizations with varying number of elements, uneven splits of the domain, as well as a discretization that contains an interior point, at which four elements intersect, (d). The wedge discretizations in Figure \ref{fig:WedgeSections} include a reference multishape, a multishape that is discretized in the radial direction as well as one in the angular direction, and a discretization containing three elements. 

The three differentiation operators, that is gradient, divergence, and Laplacian, are validated by comparing against an exact solution. The test function considered is
\begin{align*}
	g_1(x_1,x_2) &= \exp\left(-\frac{(x_1x_2 - 1)^2}{20}\right)\sin\left(\frac{x_1x_2}{4}\right),
\end{align*}
which is also used to test the accuracy of the interpolation matrix. This is tested by interpolating $g_1(x_1,x_2)$ onto a (uniform) target grid defined by $x_{1,U}$, $x_{2,U}$ and comparing the result to $g_1$ evaluated directly on the target grid. The choice of $g_1$ is largely arbitrary; similar results were obtained for a
wide range of functions.  This is also the case with the functions chosen below for other validation tests.
For integration we choose the (Gaussian) test function
\begin{align*}
	g_2(x_1,x_2) &= \exp\left(-x_1^2 - x_2^2\right).
\end{align*}
To test the convolution operator we have two different test functions
\begin{align*}
	\chi_C(\vec x) &= \exp\left(\frac{x_1^2 - x_2^2}{10}\right), \quad\hspace{0.3em} n_C(\vec z) = z_1^2 + z_1z_2,\\
	\chi_P(\vec x) &= \exp(x_1 + x_2), \ \  \qquad n_P(\vec z) = \exp\left(-z_1^2 - z_2^2 + z_1 + z_2\right),
\end{align*}
where $\chi_C$, $n_C$ are considered for Cartesian discretizations and $\chi_P$, $n_P$ are the test functions on polar multishapes.  This is 
to facilitate the analytical computation of the convolutions.

\begin{figure}[h]
	\centering
	\includegraphics[scale=0.1]{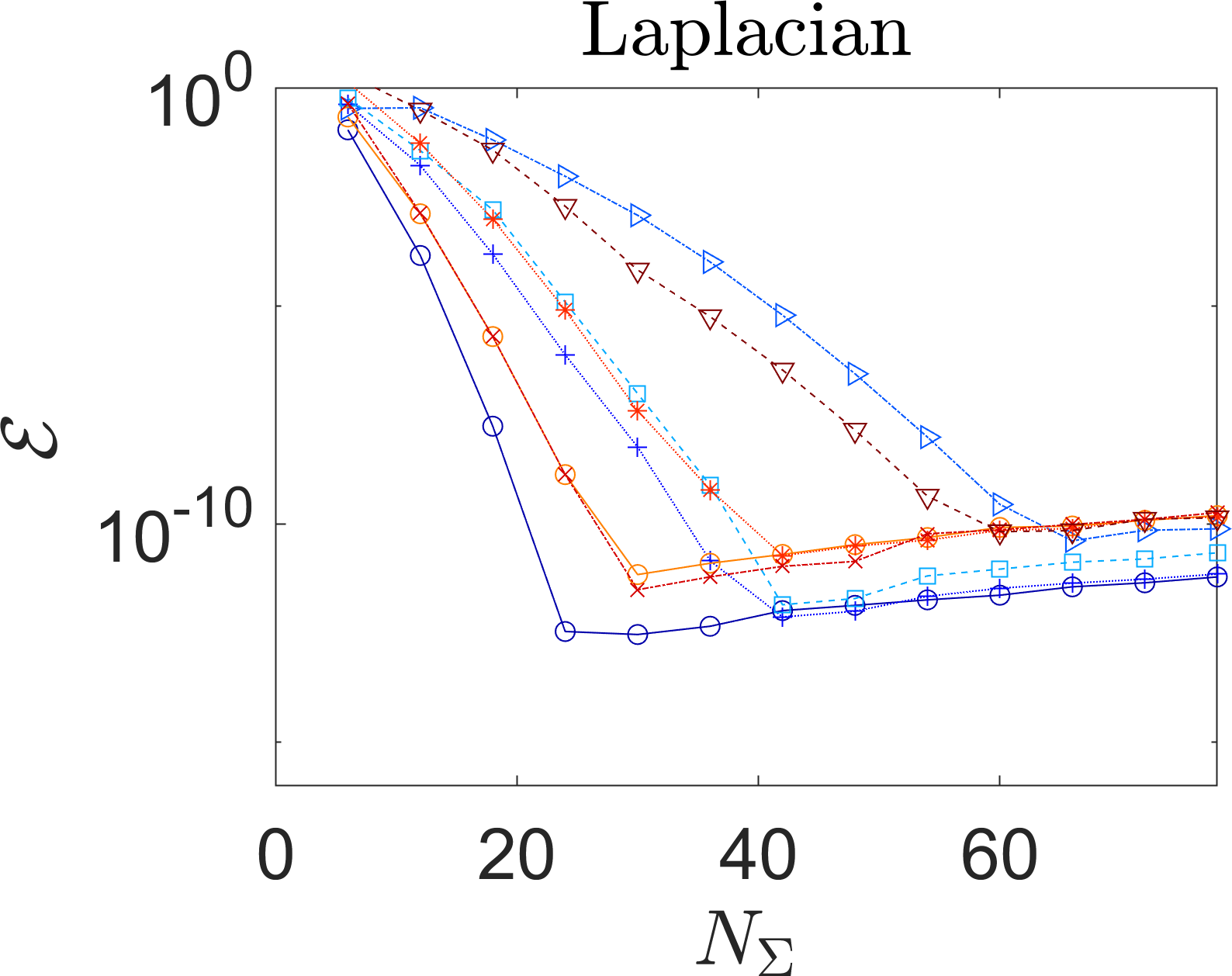}
	\includegraphics[scale=0.1]{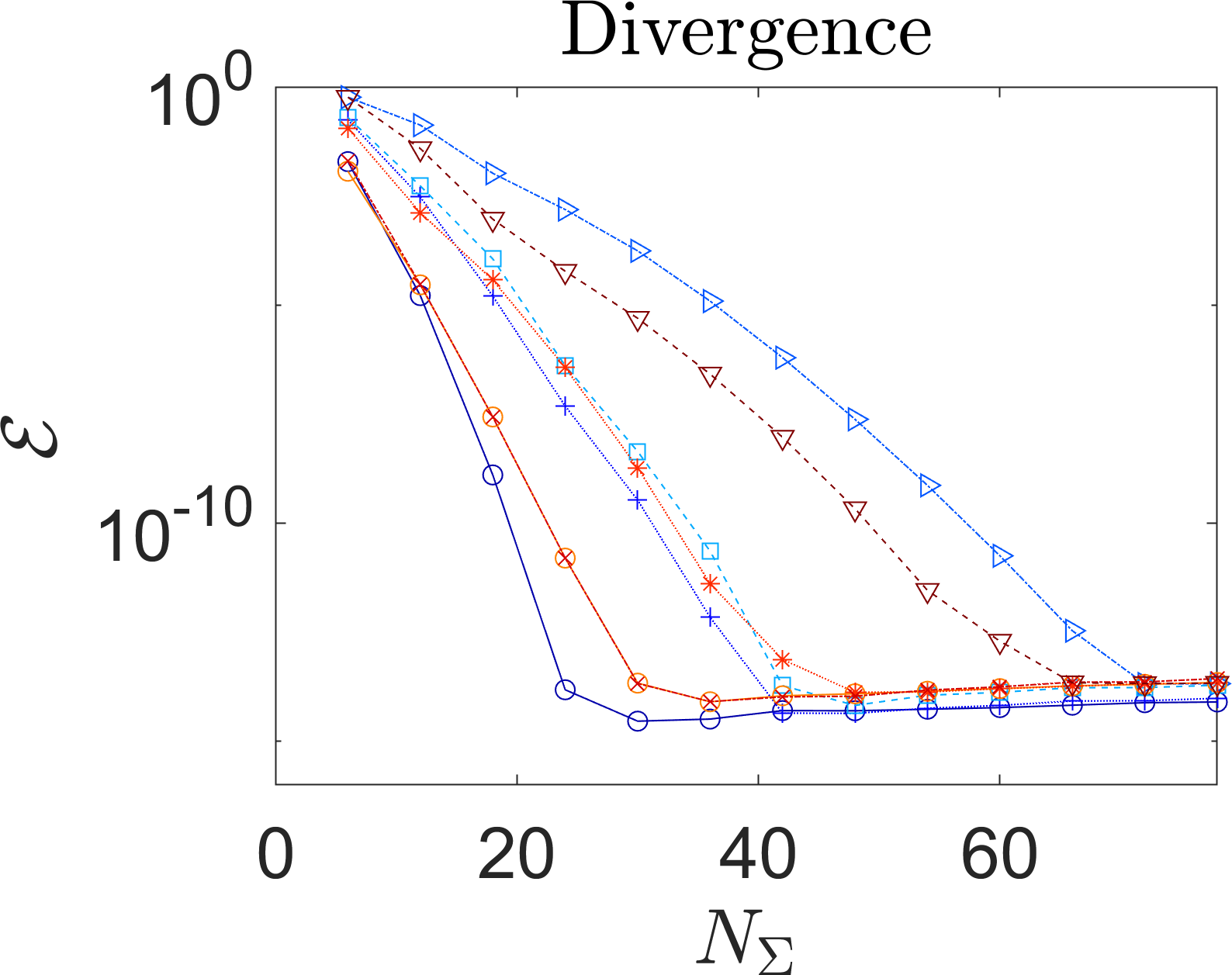}
	\raisebox{1.4cm}{\includegraphics[scale=0.15]{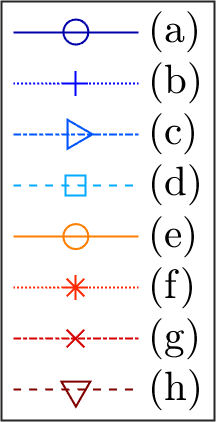}}\\
	\includegraphics[scale=0.1]{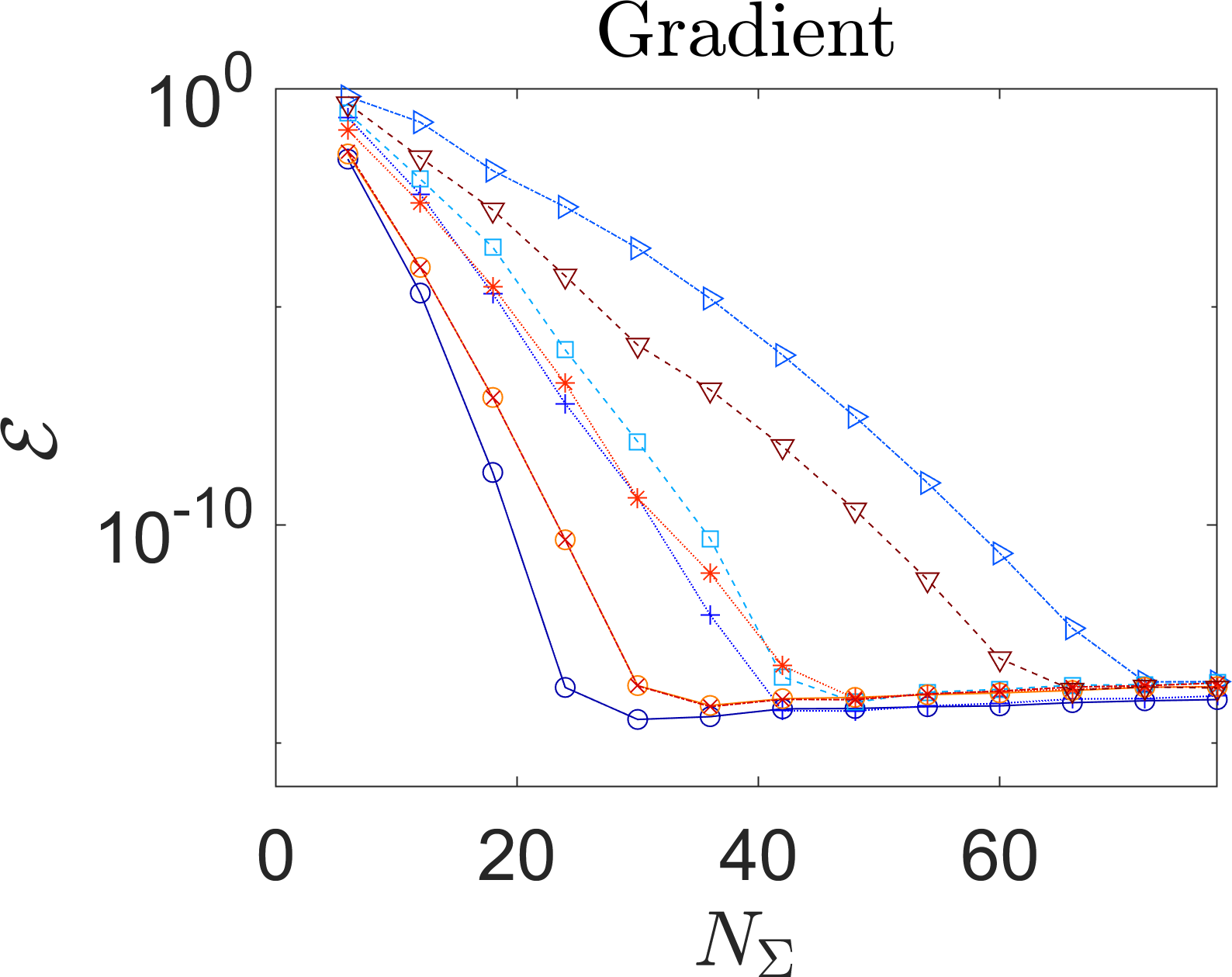}
	\includegraphics[scale=0.1]{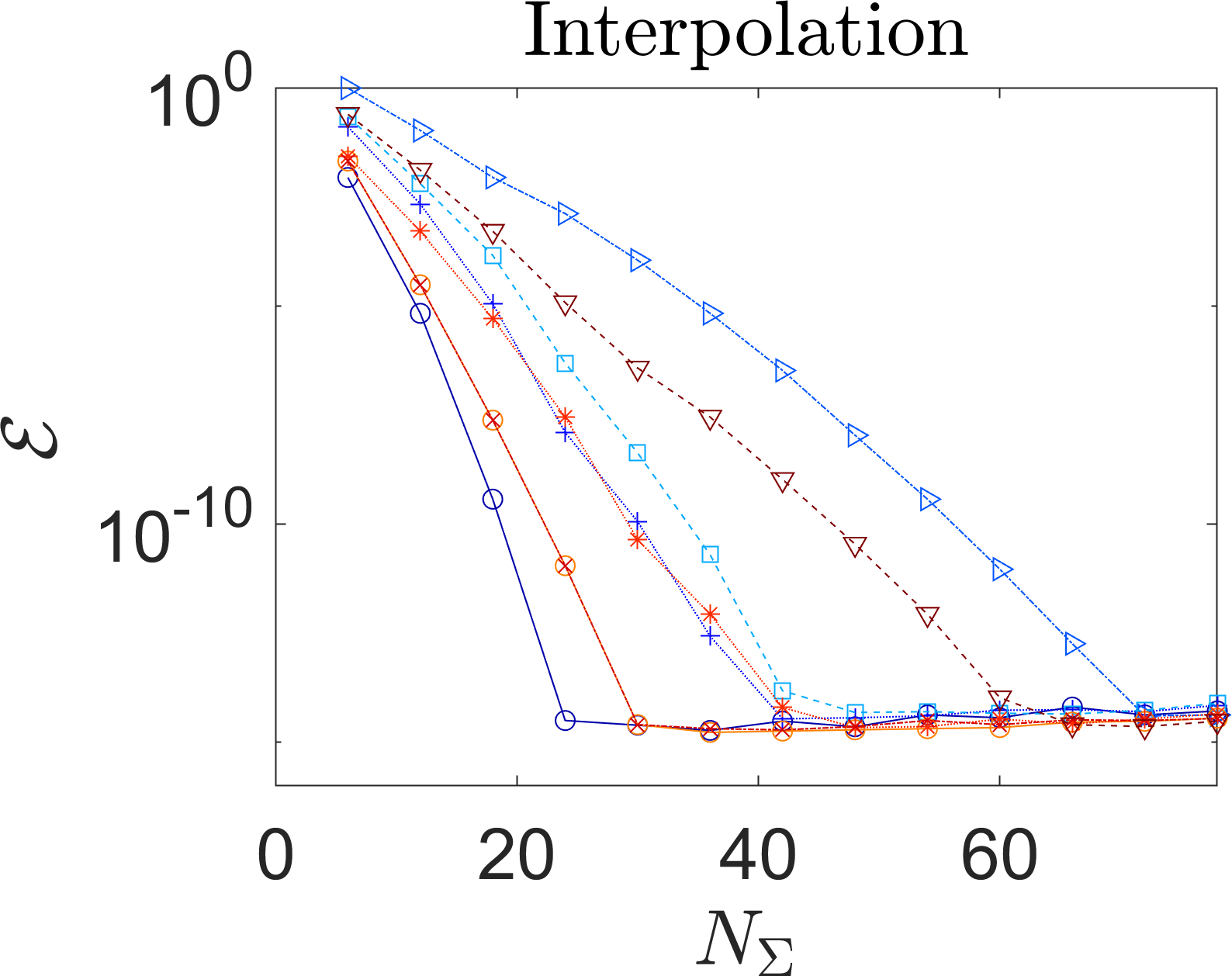}
	\phantom{\includegraphics[scale=0.15]{Legend.png}}\\
	\includegraphics[scale=0.1]{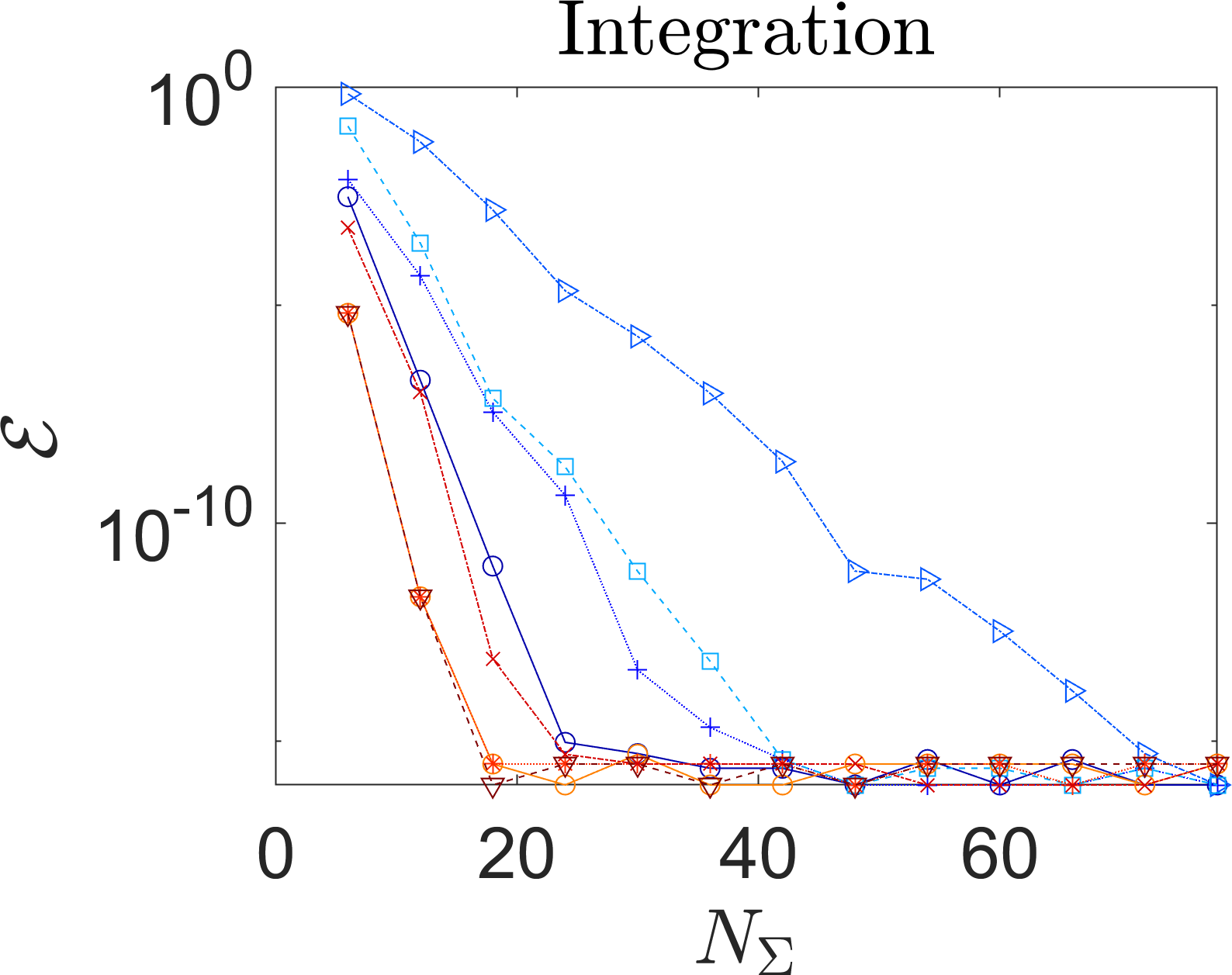}
	\includegraphics[scale=0.1]{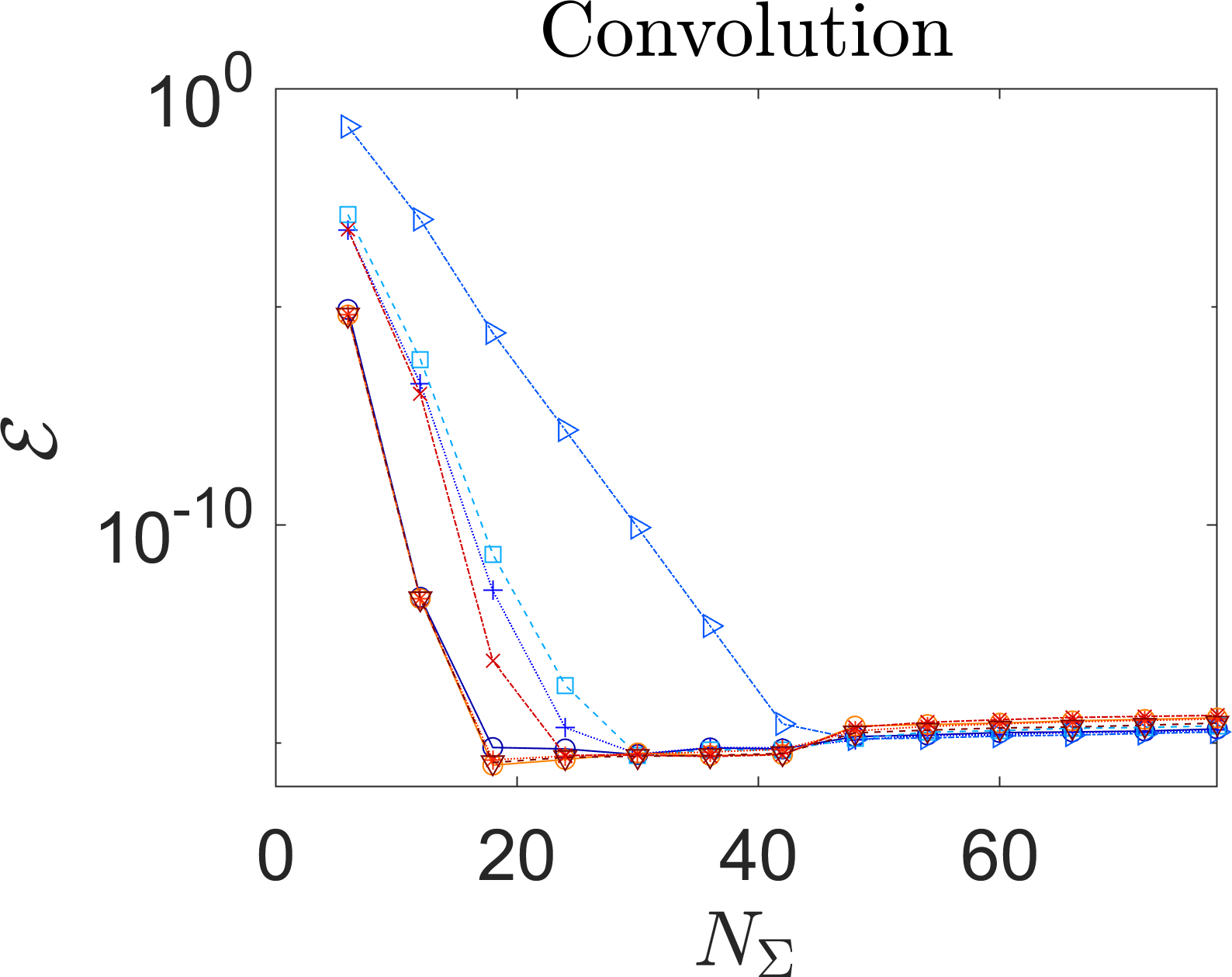}
	\phantom{\includegraphics[scale=0.15]{Legend.png}}
	\caption{The error, as defined by \eqref{eq:ErrorMeasureL2}, in the numerical differentiation, integration, interpolation, and convolution operators. Errors are computed for varying numbers of discretization points $N_\Sigma$ and measured against an exact solution.} 
	\label{fig:DiffTests}
\end{figure}

The convergence resulting from computations with the Laplacian, divergence, gradient, and interpolation matrices on different discretizations can be seen in Figure \ref{fig:DiffTests}. As expected, convergence to the exact solution is observed as the number of discretization points is increased. For each discretization, the total number of points in each of the $x_1$ and $x_2$ directions (over the whole multishape), denoted by $N_\Sigma$, remains the same. Since the overall number of points for each multishape is the same, the computational cost for each discretization is fixed and so the performance of each discretization choice can be compared. For each multishape, the $N_\Sigma$ points are equally distributed among the different elements. However, as the following comparison demonstrates, a minimum number of points has to be allocated in each the $x_1$ and $x_2$ directions and on each element of the multishape, in order to achieve optimal convergence. Multishape (c) in Figure~\ref{fig:BoxSections} is discretized into three elements in such a way that the number of points in the $x_1$ coordinate is split three ways, resulting in each element only having $N_\Sigma/3$ points in the $x_1$ coordinate, and subsequently a poorer overall convergence than, for example, in multishape (d). While multishape (d) is discretized into four shapes, the $x_1$ and $x_2$ coordinate are only divided once, so that each element has $N_\Sigma/2$ points in each direction, improving convergence considerably.

\subsection{Poisson Equation}
\begin{figure}[h]
	\centering
	\includegraphics[scale=0.13]{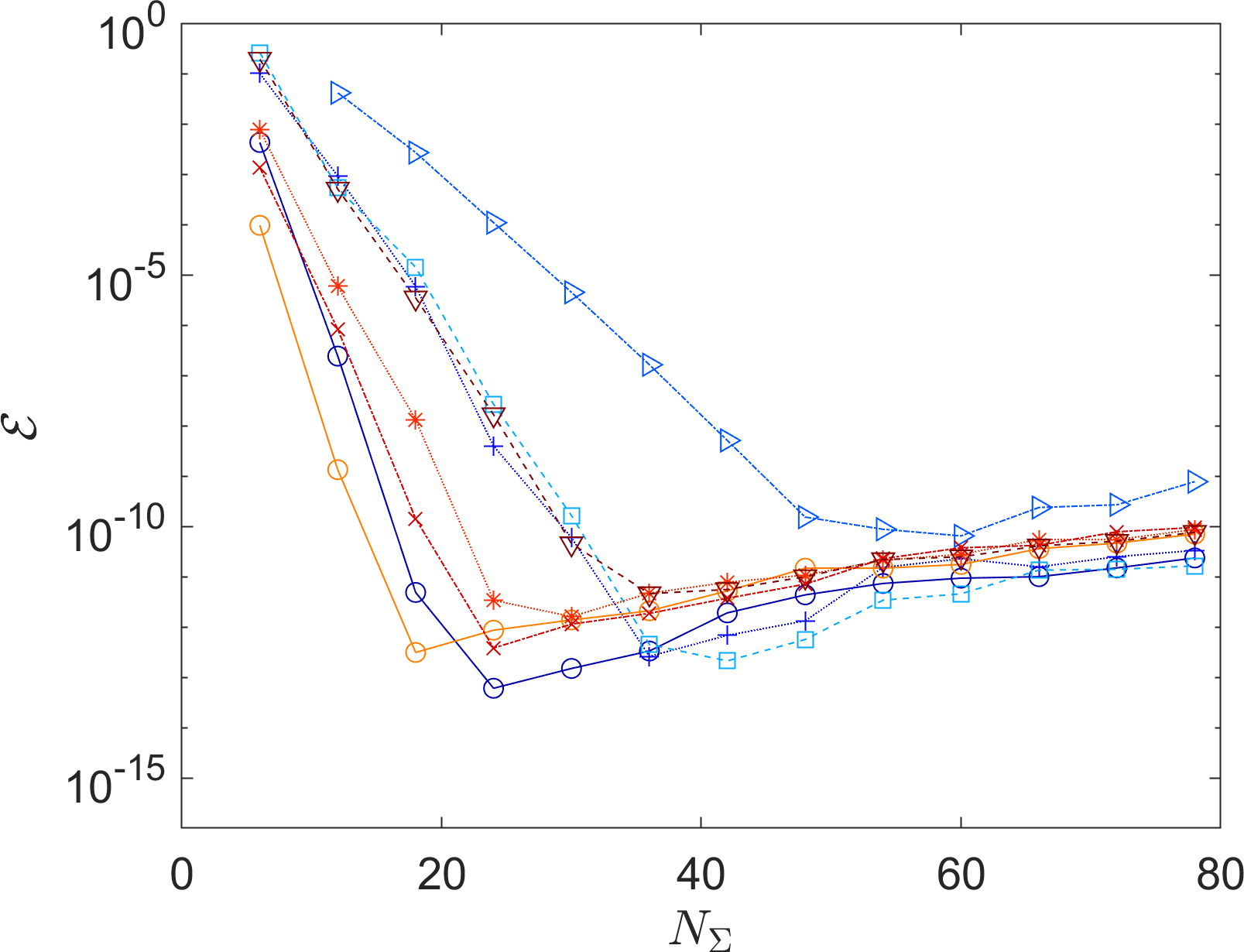}
	\raisebox{1.9cm}{\includegraphics[scale=0.18]{Legend.png}}
	\caption{The error, as defined by \eqref{eq:ErrorMeasureL2}, in computing the solution to the Poisson equation on different discretizations of a box and a wedge. Errors are computed for varying numbers of discretization points $N_\Sigma$ and measured against an exact solution. [Note that for (c), with $N_\Sigma = 6$ no interior points in $x_1$ are present, so this data point is omitted.]} 
	\label{fig:PoissonTestBW}
\end{figure}
\begin{figure}[h]
	\centering
	\includegraphics[scale=0.09]{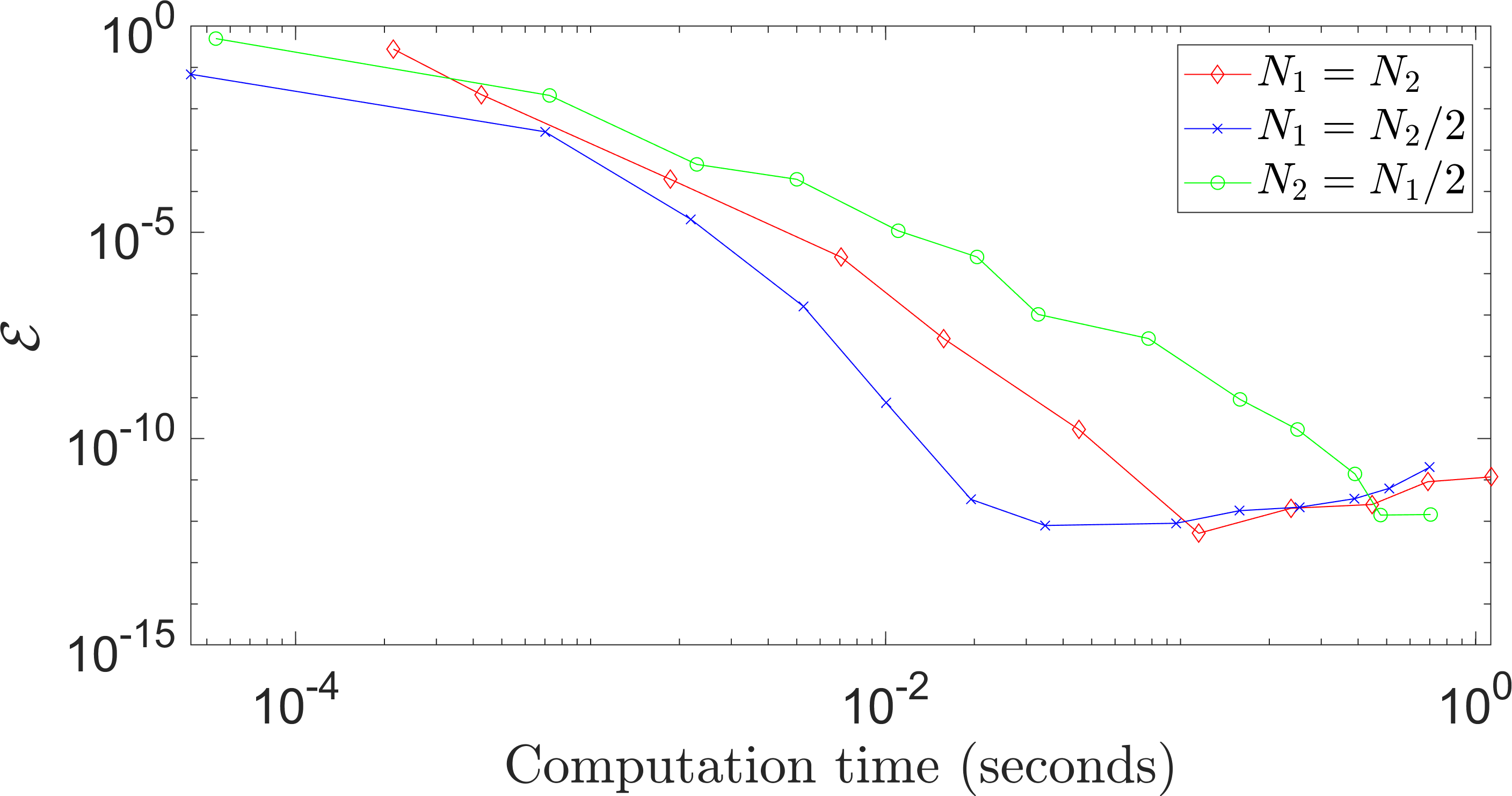}
	\caption{The error, defined by \eqref{eq:ErrorMeasureL2}, in solving the Poisson equation on a multishape, against computation time, for varying numbers of discretization points in the two coordinate directions, $N_1$ and $N_2$, for each element.} 
	\label{fig:PoissonTest}
\end{figure}
Having validated the convergence of different operators on different multishapes, for the next test we solve a Poisson equation. This demonstrates that applying the intersection matching conditions between the elements of a multishape works as expected, without adding further complexities, such as PDE solvers or integral terms.
We solve the Poisson equation 
\begin{align}\label{eq:PoissonProblem1}
	\nabla^2 \rho(\vec x) &= f(\vec x) \qquad\text{in }\Omega,\\
	\rho(\vec x) &= u(\vec x) \qquad\text{on }\partial \Omega. \notag
\end{align}
The right-hand side of the PDE,
\begin{align*}
f(\vec x) = \left(\left(x_1^2 - x_1 - 0.75\right) + \left(x_2^2 - x_2 -0.75\right)\right)\exp\left(-0.5(x_1-0.5)^2 - 0.5(x_2-0.5)^2\right),
\end{align*}
is chosen such that the Poisson equation on an infinite domain
\begin{align*}
	\nabla^2 u(\vec x) &= f(\vec x) \qquad\text{in }\Omega,\\
	u(\vec x) &\to 0 \qquad \hspace{1.225em} \text{as } |\vec x| \to \infty.
\end{align*}
has an exact solution, $u(\vec x) = \exp(-0.5(x_1-0.5)^2 - 0.5(x_2-0.5)^2)$.  This can then be used to test the accuracy for the numerical solution to the Poisson problem \eqref{eq:PoissonProblem1}, by enforcing the non-constant Dirichlet boundary conditions $\rho(\vec x) = u(\vec x)$ on the multishape boundary.  
As in Section \ref{sec:ValidationOP}, we test the convergence on the different discretizations of the box (Figure \ref{fig:BoxSections}) and wedge (Figure \ref{fig:WedgeSections}), for varying number of discretization points, as displayed in Figure \ref{fig:PoissonTestBW}. Similarly, the solution to the Poisson equation on a multishape, which combines a quadrilateral and a wedge element, as in Figure \ref{fig:MS1}, is computed. 
Intersection matching conditions are applied between the two elements as illustrated in Section \ref{sec:ClassMethods}. Note that the two elements are not defined in the same coordinate system, which MultiShape automatically takes into account when applying the intersection condition. The three test cases considered involve taking equal numbers of points in both coordinate directions, half the number of points in the first coordinate direction ($x_1$ and $r$, respectively), and half the number of points in the second coordinate direction ($x_2$ and $\theta$, respectively). We compare the error in the numerical solution to the Poisson problem \eqref{eq:PoissonProblem1} with the resulting computational timings. Figure \ref{fig:PoissonTest} showcases the convergence of the numerical result to the exact solution with increasing time, which correlates with increasing total numbers of discretization points. The results in Figure \ref{fig:PoissonTest} demonstrate that the distribution of points among the coordinate directions does not result in a large difference in computation time taken for comparable accuracy in the resulting solution. However, the most efficient discretization for the particular choice of multishape and Poisson problem appears to be distributing double the points in the $x_2$ and angular directions than in the $x_1$ and radial directions.

\section{Applications and Model Equations}
\label{sec:ModelEquations}
An important motivation for developing this MultiShape methodology is its application to Dynamic Density Functional Theory and optimal control problems.
For simple shapes, using 2DChebClass, some such applications were tackled successfully in previous work \cite{GoddardPseudospectralCode1,OurPreprint2021}. In the following subsections, we will introduce the necessary theory, relevant models, and additional numerical methods used in combination with the MultiShape implementation. 

\subsection{A Dynamic Density Functional Theory Model} \label{sec:DDFTModel}

Dynamic Density Functional Theory (DDFT) describes how a one-body particle density $\rho$ changes depending on time and the free energy of the system $\mathcal F$. This results in, often complicated, integro-PDE models. A good overview of the various applications of the theory is provided in \cite{DDFTReview}. In this work we focus on a DDFT model that includes particle interactions via a mean-field model, and leave more complex formulations that describe effects such as (hard-core) volume exclusion for future study. 
This idea is easily extended to a system with more than one particle species.
The  $n_s$  integro-PDEs that describe the dynamics of a system with $n_s$ species of particles, with densities $\rho_a$, for $a = 1,...,n_s$, 
are of gradient flow form
\begin{align} \label{eq:DDFTProblem1}
	\frac{\partial \rho_a}{\partial t} = \nabla \cdot \left( \rho_a \nabla \frac{\delta \mathcal F(\{\rho_b\})}{\delta \rho_a}\right),
\end{align}
where $\mathcal{F}$ is the Helmholtz free energy functional depending on all species $b = 1,..,n_s$.\footnote{For consistency with the cost functional \eqref{eq:OCPSys}, we adopt the notation $\mathcal{F}(\{\rho_b\})$. However, we note that $\mathcal{F}[\{\rho_b\}]$ is the standard notation in the DDFT community, explicitly denoting the \emph{functional} dependence on $\{\rho_b\}$. } The free energy of a system is composed of an ideal gas contribution, as well as contributions from an external potential and from particle--particle interactions. Under the mean-field approximation, such a system has a free energy functional given by
\begin{align*}
	\mathcal F(\{\rho_b\}) =& \sum_{a=1}^{n_s} \int_\Omega \rho_a (\ln \rho_a - 1) d \vec x + \int_\Omega \rho_a V_{\text{ext}} d\vec x\\
	&+ \frac{1}{2} \sum_{a=1}^{n_s}\sum_{b=1}^{n_s} \int_\Omega \int_\Omega \rho_a(\vec x) \rho_b(\vec x') V_{a,b}(|\vec x - \vec x'|) d\vec x' d\vec x, \notag
\end{align*}
where $V_{\text{ext}}$ is an external potential and $V_{a,b}$ a potential, describing the interactions between species $a$ and $b$. Eq. \eqref{eq:DDFTProblem1} is exactly \eqref{eq:ProblemStatement1} introduced in Section \ref{sec:ModelIntro}, with $- \mathfrak j_a = \rho_a \nabla   \frac{\delta \mathcal F}{\delta \rho_a}$. In this work, the interaction potential is given by 
\begin{align}\label{eq:V2}
	V_{a,b}(x) = \kappa_{a,b}\exp\left(- (x/\sigma_{a,b})^2\right),
\end{align}
where $\sigma_{a,b}$ represents half of the distance between the centres of particle $a$ and particle $b$ and with interaction strength $\kappa_{a,b}$. If $a = b$, this can be thought of as the particle radius. Again, the results of the method are highly robust under
different choices of $V_{a.b}$ (not shown here).  Note that here the de Broglie wavelength $\Lambda = 1$, and $k_BT = 1$.
After inserting the free energy into \eqref{eq:DDFTProblem1} we get a system of PDEs of the form
\begin{equation}\label{eq:FWPDEn}
	\frac{\partial \rho_a}{\partial t} = \nabla^2 \rho_a + \nabla \cdot \left( \rho_a \nabla V_{\text{ext}}\right) + \sum_{b=1}^{n_s} \nabla_{\vec x} \cdot \int_\Omega \rho_a(\vec x) \rho_b(\vec x') \nabla_{\vec x'} V_{a,b} (|\vec x-\vec x'|)d\vec x',
\end{equation}
with appropriate initial and boundary conditions. We focus here on no-flux boundary conditions, since they are some of the most complicated conditions to apply, in particular for non-local problems, and are highly relevant for applications. These are defined as
\begin{align}\label{eq:nofluxn}
	&  \rho_a \nabla \frac{\delta \mathcal F}{\delta \rho_a} \cdot \vec n =  \frac{\partial \rho_a}{\partial n} + \rho_a \frac{\partial V_{\text{ext}}}{\partial n} +\sum_{b=1}^{n_s}  \int_\Omega \rho_a(\vec x) \rho_b(\vec x') \frac{\partial V_{a,b}}{\partial n}(|\vec x-\vec x'|) d\vec x'  = 0\notag
\end{align}
on $\partial \Omega$, where $\vec n$ is the outward normal and $\frac{\partial}{\partial n}$ denotes the partial derivative with respect to this normal. To add some additional complexity, we consider the particle dynamics of two species; the generalisation to more species is straightforward. This system of integro-PDEs will be solved using {\scshape Matlab}'s inbuilt DAE solver \texttt{ode15s}, which is demonstrated in detail in Listing \ref{code:PDEsol}. An illustrative numerical experiment for a DDFT problem is shown in Section \ref{sec:DDFTResult}.

\subsection{Equilibrium Solution of a DDFT Model} \label{sec:EqDFTModel}
Often one is interested not (only) in the time-dependent solution to the DDFT equations, but (also) in the long-time, equilibrium solution. We would like to find such a solution on multishape domains. There are two equivalent approaches to finding these steady states; investigating the behaviour of \eqref{eq:DDFTProblem1} for large $t$, such that $\pd{\rho_a}{t} \approx 0$, or minimizing $\mathcal F$ with respect to $\rho_a$. 
Focusing on the latter approach, we solve the system of equations
$\frac{\delta \mathcal F}{\delta \rho_a} = 0$,
resulting in the non-linear equation
\begin{align}\label{eqn:eqDFT1}
	\rho_a = Z^{-1} \exp \left(-V_{\text{ext}} - \sum_{b=1}^{n_s} \int_\Omega \rho_b(\vec x')V_{a,b}(|\vec x-\vec x'|)d\vec x'\right),
\end{align}
 where $Z^{-1}$ is a normalization constant, enforcing that the mass of, and therefore the number of particles in, the system is fixed. This system of equations is of the form $\rho_a = G_a(\rho_1,...,\rho_{n_s})$ and can be solved self-consistently using, for example, iterative methods. As above, in this work, we are considering a system of two species, which demonstrates all relevant features of the numerical solution, but this can can easily be extended to systems with more species. In order to find the solution to the self-consistent system of equations \eqref{eqn:eqDFT1}, we use an iterative scheme, known as Picard iteration in the DDFT community \cite{roth2010fundamental}. Algorithm \ref{alg:PicardEqSolver} illustrates this scheme, given initial guesses for the particle distributions $\rho_{1,\text{ig}}$, $\rho_{2,\text{ig}}$, an error tolerance to terminate the algorithm, and a mixing parameter $\lambda$ for the Picard step. The Picard step takes a linear combination of the previous and new guesses for $\rho_1$, $\rho_2$, which is designed to stabilize the algorithm. 
The error measure in this algorithm is defined as the maximum of the error across all species, individually computed using \eqref{eq:ErrorMeasureL2}. 
A numerical example of this formulation can be seen in Section \ref{sec:DFTResult}.
\begin{algorithm}[H]
	\begin{algorithmic}[1]
		\STATE Set $\rho_{1,\text{ig}}$, $\rho_{2,\text{ig}}$, Tol, and $\lambda$.
		\STATE Set $\rho_{1,\text{old}} = \rho_{1,\text{ig}}$, $\rho_{2,\text{old}} = \rho_{2,\text{ig}}$.
		\WHILE{$\mathcal E < $ Tol}
		\STATE Compute $\rho_{1,\text{new}} = G_1(\rho_{1,\text{old}},\rho_{2,\text{old}})$, $\rho_{2,\text{new}} = G_2(\rho_{1,\text{old}},\rho_{2,\text{old}})$.
		\STATE Compute error $\mathcal E$, defined in \eqref{eq:ErrorMeasureL2}.
		\IF {$\mathcal E \geq $ Tol}
		\STATE Set $\rho_{1,\text{update}} = (1- \lambda)\rho_{1,\text{old}} +\lambda\rho_{1,\text{new}}$, $\rho_{2,\text{update}} = (1- \lambda)\rho_{2,\text{old}} +\lambda\rho_{2,\text{new}}$.
		 \STATE Update $\rho_{1,\text{old}} = \rho_{1,\text{update}}$, $\rho_{2,\text{old}} = \rho_{2,\text{update}}$.
		 \ENDIF
		\ENDWHILE	
	\end{algorithmic}
	\caption{Picard Equilibrium Solver}
	\label{alg:PicardEqSolver}
\end{algorithm} 
Note that problems in equilibrium Density Functional Theory (DFT) also require finding a minimum of a free energy $\mathcal F$, and can therefore be tackled by such variational approaches. An implementation is discussed in \cite{roth2010fundamental}, while overviews of DFT are given in, e.g., \cite{evans1979nature,lowen1994melting,wu2006density,wu2007density}.

\subsection{Optimal Control with Non-Local PDE Constraints}\label{sec:OCPModel}
Finally, we are interested in solving an optimal control problem involving a system of two interacting species of particles, such as in \eqref{eq:FWPDEn}. Here the challenge is to apply a control to the system, which is governed by the set of PDEs, so that the state approximates a desired `target' state, whilst also balancing the `cost' of the control. This leads to a particularly complicated system of four PDEs, with several convolution integrals describing the different interactions. This can be solved by combining MultiShape with a recently-developed Newton--Krylov algorithm \cite{guttel2021spectral} for non-linear optimization problems. The extended formalism to tackle non-local problems has been introduced in \cite{OurPreprint2021} for single species systems on rectangular domains and is extended in this section to multiple species on a multishape domain. Overviews of PDE-constrained optimization are given in, e.g., \cite{hinze1,delosreyes1,troeltsch1}.
For a system of $n_s$ species, we can formulate the following optimal control problem
\begin{equation}\label{eq:OCPSys}
	\min_{\{\rho_a\}, \vec w} \mathcal J (\{\rho_a\}, \vec w) := \frac{1}{2}\sum_{a=1}^{n_s}\int_0^T \int_\Omega \left(\rho_a - \hr_a \right)^2 d\vec x dt + \frac{\beta}{2} \int_0^T \int_\Omega |\w|^2 d\vec x dt,
\end{equation}
subject to the following system of $n_s$ PDE constraints, for  $a = 1,...,n_s$:
\begin{align*}
	\frac{\partial \rho_a}{\partial t} = \nabla^2 \rho_a - \nabla \cdot \left( \rho_a \vec w\right) + \nabla \cdot \left( \rho_a \nabla V_{\text{ext}}\right) \qquad\qquad\qquad\qquad\qquad& \\
	+ \sum_{b=1}^{n_s} \nabla_{\vec x} \cdot \int_\Omega \rho_a(\vec x) \rho_b(\vec x') \nabla_{\vec x'} V_{a,b} (|\vec x -\vec x'|)d\vec x' \quad &\text{in } \Omega \times (0,T), \\
	\frac{\partial \rho_a}{\partial n} - \rho_a \vec w \cdot \vec n + \rho_a \frac{\partial V_{\text{ext}}}{\partial n}  + \sum_{b=1}^{n_s}  \int_\Omega \rho_a(\vec x) \rho_b(\vec x') \frac{\partial V_{a,b}}{\partial n}(|\vec x-\vec x'|) d\vec x'  = 0 \quad &\text{on } \partial \Omega \times (0,T),\\
	\rho_a(\vec x, 0) = \rho_{a,0}(\vec x) \quad&\text{at } t = 0,
\end{align*}
where $T>0$.
Note that the difference between the system of PDE constraints above and the system \eqref{eq:FWPDEn} is the additional term involving $\vec w$. This term models how a particle distribution $\rho_a$ may be influenced by a background flow $\vec w$. In this optimal control problem, the aim is to drive the particle distributions $\rho_a$ closer to some desired states $\hr_a$, which can be achieved by modifying $\vec w$, the control variable.
Note that in DDFT, $\vec w$ is generally of the form $\vec w = \nabla V$, where $V$ is an external potential, see \cite{DDFTReview} for a discussion on considerations of other forms of $\vec w$ and their limitations. Enforcing this condition in an optimal control setting is more complicated, however, since enforcing the gradient of the quantity we aim to control is a problem without a unique solution. At present, we restrict ourselves to control via a general $\vec w$ as a test case for the methodology.
The regularization parameter $\beta$ describes how much the use of control is penalized when solving the system and consequently how close each $\rho_a$ can be driven to the targets $\hr_a$. This minimization is subject to the physical constraints described by the system of PDEs and the corresponding boundary conditions. 
Such an optimization problem can be solved by considering the corresponding first-order optimality system, applying the working in \cite[Section 3]{OurPreprint2021}. It consists of the system of PDE constraints and the adjoint PDE system where, for each species $a$, we have the following PDE for the adjoint variable $\adj_a$:
\begin{align*}
	\frac{\partial \adj_a}{\partial t} =& - \nabla^2 \adj_a - \w \cdot \nabla \adj_a + \nabla V_{\text{ext}} \cdot \nabla \adj_a  + \hr_a - \rho_a \\
	&+\sum_{b=1}^{n_s}  \int_\Omega \rho_b(\vec x')\left( \nabla \adj_a(\vec x) - \nabla \adj_b(\vec x') \right) \cdot \nabla V_{a,b}(|\vec x-\vec x'|)d\vec x',
\end{align*}
with boundary condition $\frac{\partial \adj_a}{\partial n} = 0$ and final time condition $\adj_a(\vec x,T) = 0$. The final part of the first-order optimality system is a single gradient equation
\begin{align*}
	\w = -\frac{1}{\beta} \sum_{a=1}^{n_s}\rho_a \nabla \adj_a.
\end{align*}
This system of PDEs for the densities and adjoints ($\rho_a$ and $\adj_a$) is solved by combining the MultiShape methodology with a Newton--Krylov algorithm, introduced in \cite{guttel2021spectral} and extended to non-local problems with non-linear controls in \cite{OurPreprint2021} -- we refer to these references for detailed descriptions. The idea of the algorithm is that a residual vector function is constructed based on the integral in time of the residuals of the state and adjoint equations, defined at every time step. The integrals are approximated using Chebyshev interpolation on the time grid (also using 2DChebClass), whereupon the residual vector is supplemented by the misfit between the computed and exact initial/final-time conditions. This allows us to construct the Jacobian via spatial discretization at each time step, whereupon the resulting Newton systems are approximately solved using the GMRES method \cite{gmres} in combination with a Kronecker-product based preconditioner for the Schur complement matrix. 
 We show numerical results for such a problem in Section \ref{sec:OCPResult}.

\section{Numerical Examples} \label{sec:NumericalResults}
In this section, problems involving the models introduced in Section \ref{sec:ModelEquations} are solved numerically on a multishape, with appropriate intersection conditions applied between multishape elements.
Throughout this section, initial conditions for $\rho_a$, $a = 1,2$, are of the form
\begin{align}\label{eq:ICrho1}
	\rho_{a,\text{ic}} = \frac{c_M f_{a,\text{ic}}}{\int_\Omega f_{a,\text{ic}} d \vec x}. 
\end{align}
Note that the model problems considered in this section are presented in order of increasing computational complexity instead of in the order in which the models were introduced in Section \ref{sec:ModelEquations}. All tests are carried out on Dell PowerEdge R430 running Scientific Linux 7, four Intel Xeon E5-2680 v3 2.5GHz, 30M Cache, 9.6 GT/s QPI 192 GB RAM and are tested using Matlab Versions 2020a--2021b. 

\subsection{Equilibrium Solution of a DDFT Model} \label{sec:DFTResult}
\begin{figure}[h]
	\centering
	\includegraphics[scale=0.08]{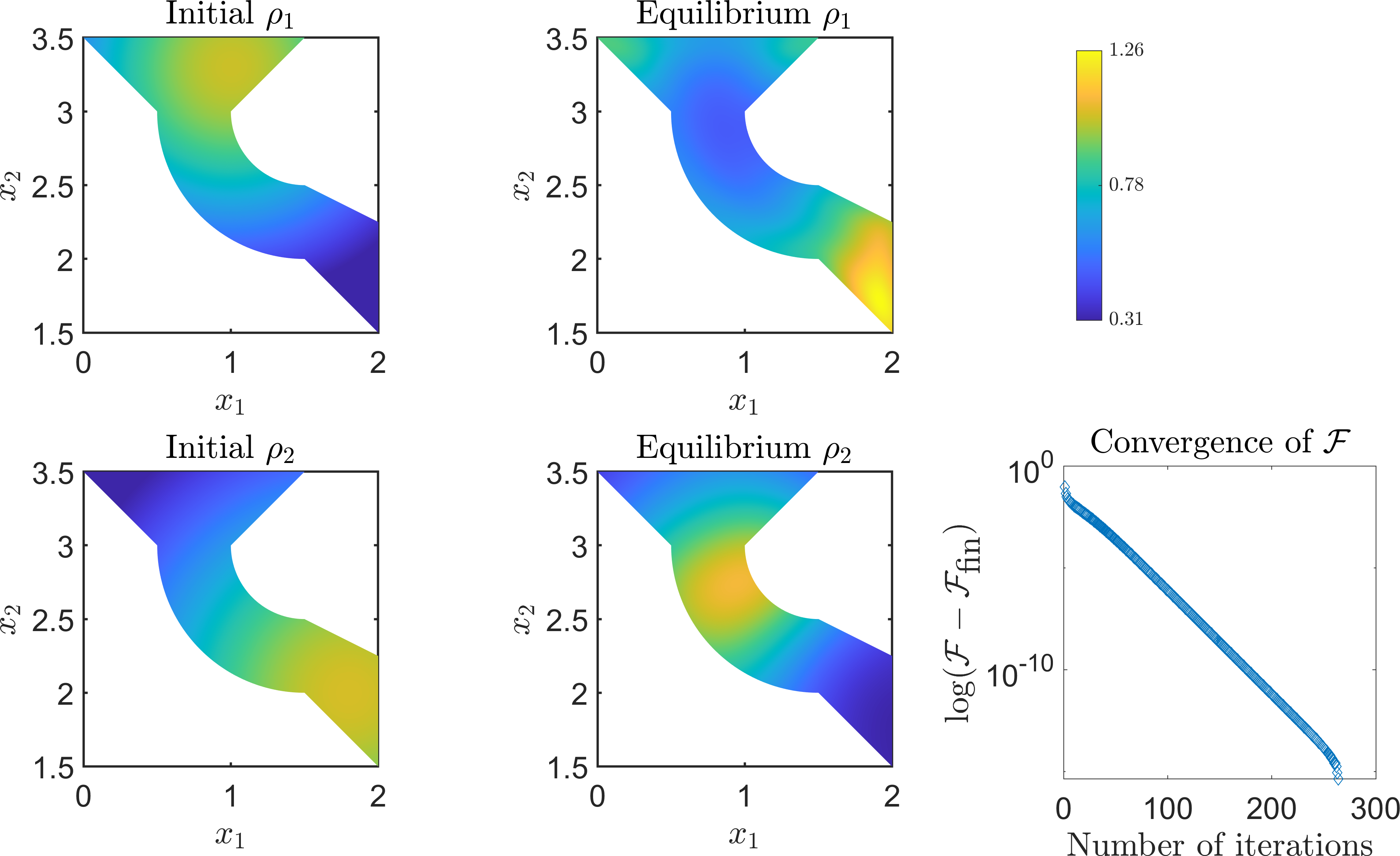} 
	\caption{Equilibrium solution for a two species system on a multishape. Bottom right: Convergence of the free energy $\mathcal F$, measuring the log change in $\mathcal F$ shifted by the final value of the free energy $\mathcal F_{\emph{fin}}$ in order to compare the differences in free energy rather than absolute values. }
	\label{fig:Eq1}
\end{figure}

We first solve an example involving the model equations introduced in Section \ref{sec:EqDFTModel}, describing the equilibrium properties of a system with two interacting particle species.
For this example, a multishape is defined using three elements: two quadrilaterals and one wedge, see Figure \ref{fig:Eq1}. The number of points on each shape is $N_1 = N_2 = 20$, chosen as a balance between accuracy and computation time, as informed by the validation tests presented in Section \ref{sec:Validation}. We choose the parameters of the problem as follows. Interaction strengths, for an interaction potential of the form \eqref{eq:V2}, are $\kappa_{1,1} = -7$, $\kappa_{2,2} = -3$, $\kappa_{1,2} = \kappa_{2,1} = 2$, with $\sigma_{1,1} = 0.1$, $\sigma_{2,2} = 1$, $\sigma_{1,2} = \sigma_{2,1} = 0.55$. This simulates Species $1$ being significantly smaller than Species $2$, with intra-species attractive potentials and inter-species repulsion. We furthermore impose an external potential $V_{\text{ext}} = 0.1x_2$, simulating gravitational forces acting on the particles. While the equilibrium results are robust for different choices of initial conditions, we choose  $\rho_{1,\text{ic}}$, $\rho_{2,\text{ic}}$ of the form \eqref{eq:ICrho1}, with $c_M = 1$ and $f_{1,\text{ic}} = \exp\left(-0.5(x_1-1)^2 - 0.5(x_2-3.3)^2\right)$, $f_{2,\text{ic}} = \exp\left(-0.3(x_1-1.8)^2 - 0.3(x_2-2)^2\right)$, since they are relatively far from both the equilibrium and each other. With this choice of initial conditions, Algorithm \ref{alg:PicardEqSolver}, with $\lambda = 0.5$ and a termination tolerance of $10^{-8}$, converges to an equilibrium state within $266$ iterations, taking $6$ seconds. The initial conditions, the equilibrium solution, as well as the convergence of the free energy of the system, are displayed in Figure \ref{fig:Eq1}. In the equilibrium solution, the interplay between the gravitational forces and particle interactions can be observed. Since Species $1$ and $2$ repel each other, they are found to separate within the domain. Furthermore, since Species $1$ is smaller in size than Species $2$, it clusters at the boundary of the domain, while Species $2$ accumulates in the middle. The effect of the gravitational potential is evident when considering the proportion of Species $1$ at the top and bottom of the domain. Since gravity is acting on the particles, it is more natural for the particles to accumulate near the bottom of the domain. Species $2$ is also subject to gravitational forces, since without gravity it would be located in the top part of the domain. Increasing the gravitational force, while keeping the interactions the same, forces Species $2$ downwards towards Species $1$, while causing a steeper accumulation of particles of Species $1$ in the bottom corner. However, this is not shown in Figure \ref{fig:Eq1} in the interest of brevity.

\subsection{Time-Dependent Solution of a DDFT Model}\label{sec:DDFTResult}
\begin{figure}[h]
	\centering
	\includegraphics[scale=0.08]{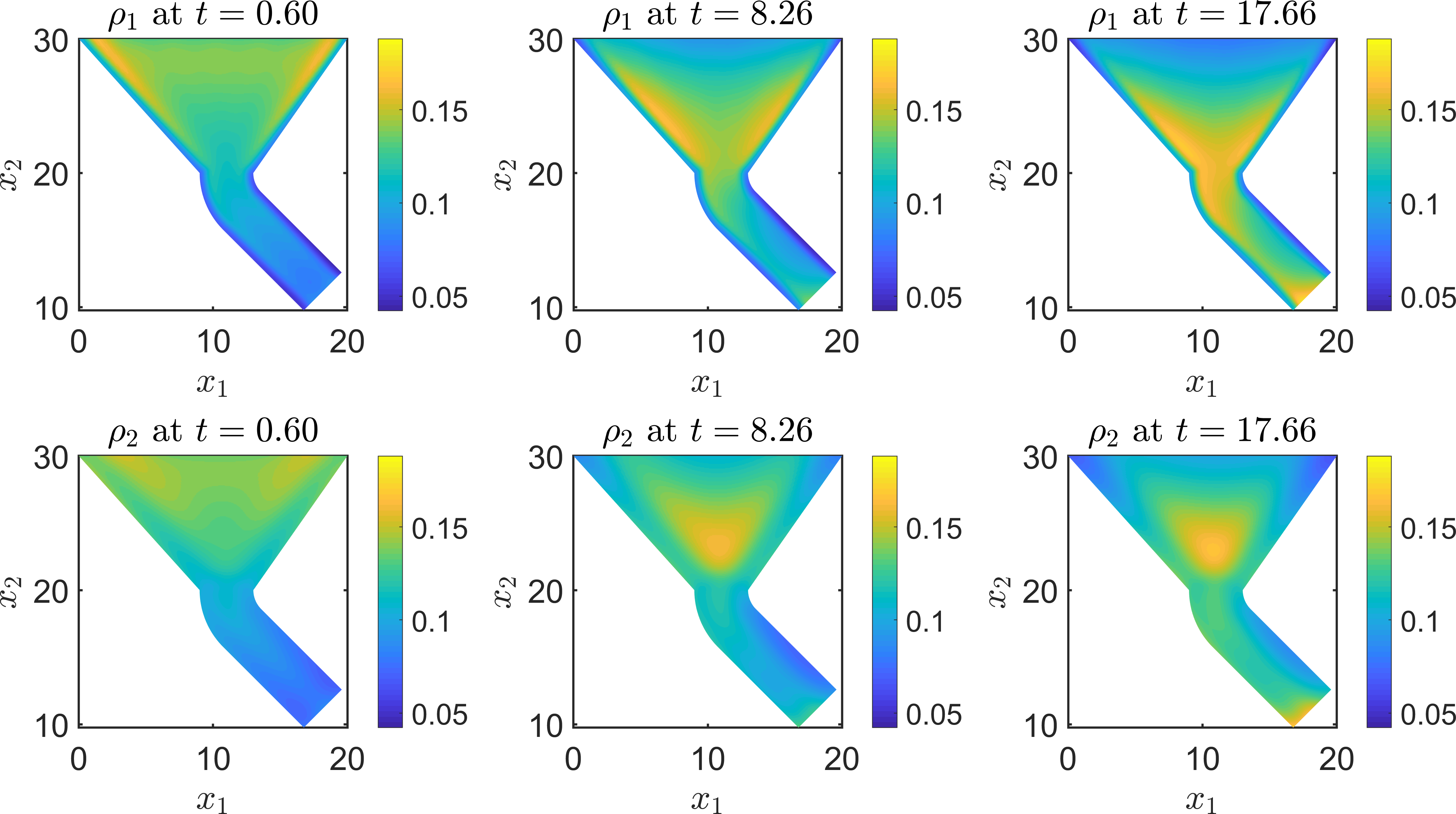} 
	\caption{Dynamics of two interacting species of different size in a funnel under the influence of gravitational forces.} 
	\label{fig:Dyn1}
\end{figure}

We now consider an example of the dynamic model introduced in Section \ref{sec:DDFTModel}, which models a funnel through which particles of different sizes are passing. There are two types of particles in the system, one species being larger than the other. We expect that smaller particles should pass through the narrow part of the funnel more easily than the larger particles. We introduce an external potential $V_g = c_g x_2$, modelling gravity, with $c_g = 0.15$.
The left and right walls of the funnel are repulsive; a simple way to approximate volume exclusion effects. Since the two species are of different sizes, the repulsion is of different strengths.  The external potential, modelling the repulsive walls, is defined as $ V_{\text{rep}} = \epsilon \left(\exp\left(- (d_L/ \alpha)^2\right) + \exp\left(- (d_R/ \alpha)^2\right)\right)$, where $d_L$ and $d_R$ are the shortest (Euclidean) distances to the left and right boundary of the multishape respectively. We choose $\epsilon = 0.6$, $\alpha_1 = 0.5$ for Species $1$, $\alpha_2 = 2$ for Species $2$. The interactions between particles are again described by the potential \eqref{eq:V2},  with $\kappa_{1,1} = \kappa_{2,2} =\kappa_{1,2} = \kappa_{2,1} = 0.1$ and corresponding $\sigma_{1,1} = 0.5$, $\sigma_{2,2} =2$, $\sigma_{1,2} = \sigma_{2,1} = 1.25$.
The choice of $\sigma_{2,2}$ is motivated by the fact that the narrow channel in the multishape has width $2$. The initial conditions for each species are of the form \eqref{eq:ICrho1}, with $f_{1,\text{ic}} = f_{2,\text{ic}} = x_1 + 5$ and $c_M = 20$. The time interval for this simulation is $t \in [0,20]$, and each element in the multishape is discretized using $N_1 = N_2 = 20$ points.
The DAE solver \texttt{ode15s}, with absolute and relative tolerances set to $10^{-9}$, takes around $20$ seconds to solve this problem and the result is displayed in Figure \ref{fig:Dyn1}.  It is evident that the behaviour of the two species is qualitatively different. Due to their size, the larger particles gather in the middle of the funnel, while the smaller particles gather at the boundaries. Both species experience the repulsion from the funnel boundary, but Species $2$ accumulates noticeably further away from the walls than Species $1$, due to their larger size. This arrangement causes the larger particles to be blocked by the smaller particles in entering the narrow part of the funnel, exacerbated by the fact that the particles repel each other. Both species are affected by gravitational forces, so that sedimentation can be observed and, over time, more particles accumulate in the bottom part of the domain. However, since the first part of the narrow channel is curved, some particles accumulate in this section of the domain before travelling downwards. 

\subsection{Optimal Control with Non-Local PDE Constraints}\label{sec:OCPResult}
\begin{figure}[h]
	\centering
	\includegraphics[scale=0.08]{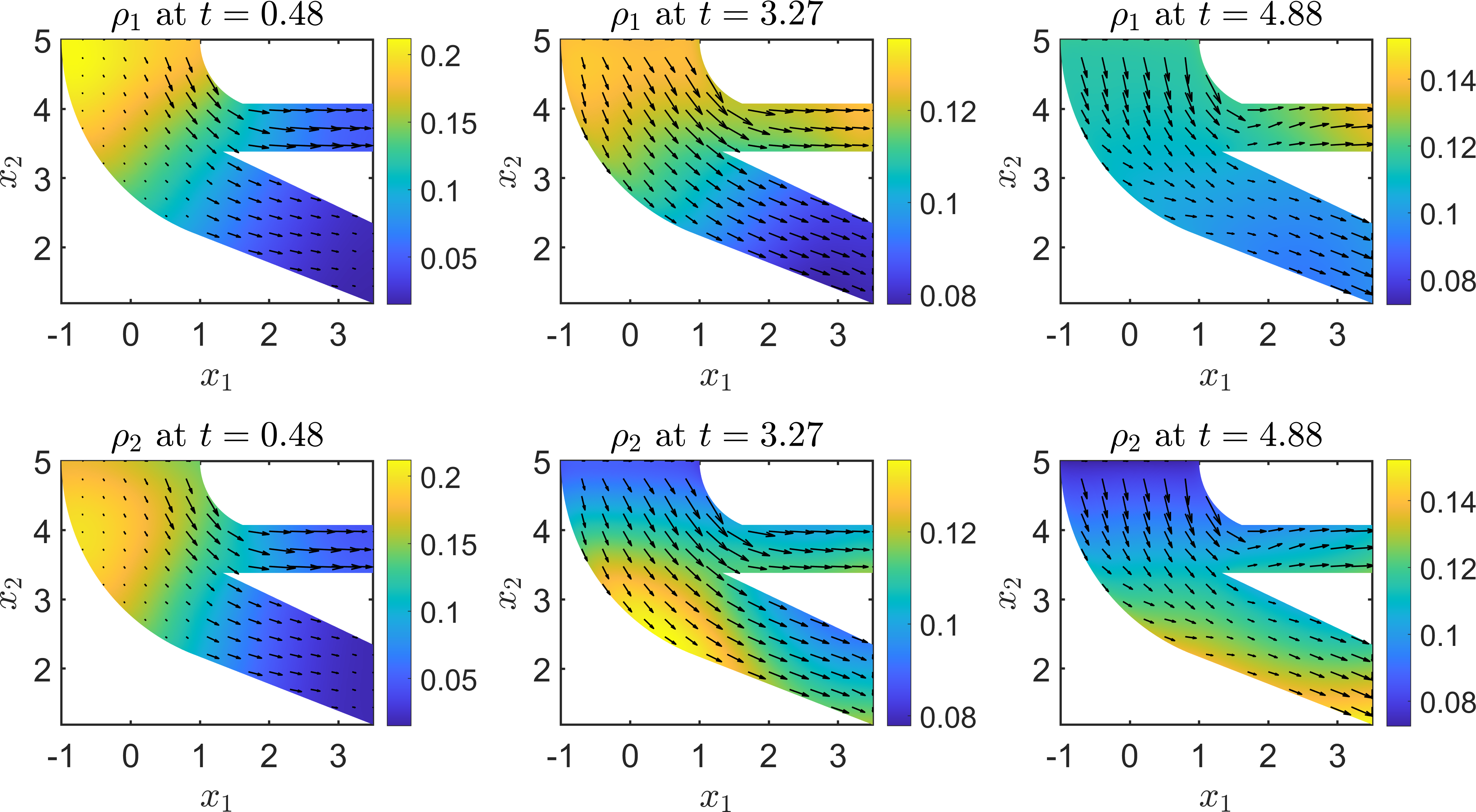} 
	\caption{Optimal states $\rho_1$ and $\rho_2$ with superimposed optimal control $\vec w$.} 
	\label{fig:OCP1}
\end{figure}

Finally, we compute results for a flow control problem of the form introduced in Section \ref{sec:OCPModel} on four elements, with $N_1 = N_2 = 14$ discretization points in space on each element and $n = 10$ points in the time interval $t \in [0,5]$. Note that the normal vector at the intersection of the four elements (see Figure \ref{fig:OCP1} for reference) has to be correctly defined in a manual fashion, as illustrated for a simpler test case in Figure \ref{fig:NormalOverride}.  We choose two interacting particle species with $\sigma_{1,1} = 0.5$, $\sigma_{2,2} = 1$, with $\kappa_{1,1} = -0.8$, $\kappa_{1,2} = 0.6$, $\kappa_{2,1} = 0.2$, $\kappa_{2,2} = -0.3$. External potentials are imposed for each species: $V_{\text{ext},1} =0.05 \left(\exp{(-0.1(x_1-2)^2 - 0.1(x_2-1)^2)} - x_2\right)$ and $V_{\text{ext},2} =0.1 \left(\exp{(-0.3(x_1-3)^2 - 0.3(x_2-4)^2)} + 2x_2\right)$. The initial conditions for $\rho_1$ and $\rho_2$ are defined by	\eqref{eq:ICrho1}, where $f_{1,\text{ic}} = f_{2,\text{ic}} = \exp{(-0.15(x_1+0.5)^2 - 0.15(x_2-5)^2)}$ and $c_M = 1$. We choose $\beta = 10^{-3}$, allowing a significant amount of control to be introduced into the system, so that we expect to get reasonably close to the target states $\hr_1$ and $\hr_2$.
The targets are defined by a forward problem with the same initial conditions and external potential, with a background flow of strength $0.1$ in the $x_1$ direction (along the channel in the wedge case). Without imposing control (i.e., with $\vec w = \vec 0$), we can evaluate the cost functional \eqref{eq:OCPSys}, which results in $\mathcal J_{\text{uc}} = 0.0038$ (`uc' referring to `uncontrolled'). The value of the cost functional for the optimized system is $\mathcal J_{\text{c}} = 1.7571 \times 10^{-4}$, showcasing that the overall cost obtained is significantly reduced when the optimal control is introduced. We can observe that both species of particles are carried along by the background flow field (optimal control) over time. Since Species $1$ experiences an upward force, it accumulates in the top-right part of the channel system at later times. Equivalently, since Species $2$ experiences a downward acting force, it accumulates in the bent part of the channel on the bottom right at later times of the simulation. It takes around $3$ hours to solve this optimal control problem, which converges within $10$ iterations of the Newton--Krylov algorithm, to a residual error of order $10^{-13}$ and $10^{-15}$ for the state and adjoint variable, respectively.
\section{Conclusions}
\label{sec:Conclusions}
In this paper, we have presented a new spectral element implementation, which allows for the efficient solution of (integro-)PDE models on complicated domains, including the application of non-linear and non-local boundary conditions. A general numerical framework has been introduced, building on the existing code library 2DChebClass for pseudospectral methods, and with a specific emphasis on user-friendly features, such as the simple application of matching and boundary conditions. After successful validation of the MultiShape methodology, various numerical examples were discussed, ranging from steady state problems, through solving (integro-)PDE models, to optimal control problems. These examples have been chosen in view of their relevance to Dynamic Density Functional Theory, which describes various industrial processes, as well as phenomena in physics and biology. There are many possible extensions to this work, such as applying the MultiShape methodology to solve more complicated DDFT models, capturing effects such as sedimentation, hydrodynamic interactions, or inertia. Further future challenges include modelling and optimization of industrial processes, for example in manufacturing and production. Other areas of research, which could potentially benefit from applying the MultiShape methodology, include models of pedestrian dynamics, brewing, and spread of diseases on complicated domains.

\section*{Acknowledgments}
JCR is supported by The Maxwell Institute Graduate School in Analysis and its Applications (EPSRC grant EP/L016508/01), the Scottish Funding Council, Heriot-Watt University, and The University of Edinburgh.
BDG acknowledges support from EPSRC grant EP/L025159/1;
RM-W was funded by EPSRC IAA PIII075;
JWP from EPSRC grant EP/S027785/1.
BDG and RM-W would like to thank Tim Hurst, Andreas Nold, and Serafim Kalliadasis for helpful discussions.

\bibliographystyle{siamplain}
\bibliography{references}
\end{document}

%% file: ex_shared.tex
\usepackage{lipsum}
\usepackage{amsfonts}
\usepackage{graphicx}
\graphicspath{{Images/}}
\usepackage{epstopdf}
\usepackage{hyperref}
\usepackage{amsmath}
\usepackage{amsfonts}
\usepackage{amssymb}
\usepackage{mathtools}
\usepackage{commath}
\usepackage{bbm}
\usepackage{algorithmic} 
\usepackage[framed,numbered,autolinebreaks,useliterate]{mcode}
\ifpdf
  \DeclareGraphicsExtensions{.eps,.pdf,.png,.jpg}
\else
  \DeclareGraphicsExtensions{.eps}
\fi

\newsiamremark{remark}{Remark}
\newsiamremark{hypothesis}{Hypothesis}
\crefname{hypothesis}{Hypothesis}{Hypotheses}
\newsiamthm{claim}{Claim}

\headers{MultiShape: A Spectral Element Method}{J. C. Roden, R. D. Mills-Williams, J. W. Pearson, and B. D. Goddard}

\title{MultiShape: A Spectral Element Method,\\
	 with applications to Dynamic Density Functional Theory and PDE-Constrained Optimization}

\author{Jonna C. Roden\thanks{The School of Mathematics and Maxwell Institute for Mathematical Sciences, The University of 	Edinburgh, Edinburgh, EH9 3FD, UK ({\tt J.C.Roden@sms.ed.ac.uk})}
	\and Rory D. Mills-Williams\thanks{University of Edinburgh, Heinrich Heine Universit\"at D\"usseldorf}
	\and John W. Pearson\thanks{The School of Mathematics and Maxwell Institute for Mathematical Sciences, The University of Edinburgh, Edinburgh, EH9 3FD, UK ({\tt j.pearson@ed.ac.uk})}
	\and Benjamin D. Goddard\thanks{The School of Mathematics and Maxwell Institute for Mathematical Sciences, The University of Edinburgh, Edinburgh, EH9 3FD, UK ({\tt b.goddard@ed.ac.uk})} }

\usepackage{amsopn}

\newcommand{\adj}{q}
\newcommand{\w}{\vec w}

\newcommand{\hr}{\widehat \rho}